\PassOptionsToPackage{colorinlistoftodos}{todonotes}
\documentclass[10pt]{article}
\usepackage[a4paper,margin=25mm]{geometry}
\usepackage{iftex}
\ifPDFTeX
  \usepackage[T1]{fontenc}
  \usepackage[utf8]{inputenc}
\fi
\usepackage{lmodern,microtype}
\usepackage{amsmath,amssymb,amsthm,mathtools,mathrsfs,bm}
\usepackage{graphicx,booktabs,array,tabularx,enumerate}
\usepackage{float,xcolor}
\usepackage{tikz}
\usetikzlibrary{arrows.meta,positioning,calc}
\usepackage[section]{placeins}
\usepackage[numbers,sort&compress]{natbib}

\usepackage[draft]{changes}
\usepackage[colorinlistoftodos]{todonotes}
\definechangesauthor[name={Daniele},color=blue]{DAN}
\setdeletedmarkup{\textcolor{red}{\sout{#1}}}
\definecolor{primary}{RGB}{20,60,110}
\usepackage[colorlinks=true,linkcolor=primary,citecolor=primary,urlcolor=primary]{hyperref}
\numberwithin{equation}{section}
\newtheorem{theorem}{Theorem}[section]
\newtheorem*{theoremA}{Theorem A}
\newtheorem*{theoremB}{Theorem B}
\newtheorem{proposition}[theorem]{Proposition}
\newtheorem{lemma}[theorem]{Lemma}
\newtheorem{corollary}[theorem]{Corollary}
\theoremstyle{definition}

\newtheorem{example}[theorem]{Example}
\theoremstyle{remark}
\newtheorem{remark}[theorem]{Remark}

\newcommand{\R}{\mathbb R}
\newcommand{\C}{\mathbb C}
\newcommand{\cH}{\mathcal H}

\newcommand{\cW}{\mathcal W}

\newcommand{\cE}{\mathcal E}

\newcommand{\Stab}{\operatorname{Stab}}
\newcommand{\sflo}{\operatorname{sf}}
\newcommand{\Mas}{\operatorname{Mas}}

\newcommand{\gr}{\operatorname{gr}}

\newcommand{\ran}{\operatorname{ran}}
\newcommand{\spec}{\operatorname{spec}}
\newcommand{\sgn}{\operatorname{sgn}}
\newcommand{\ceil}[1]{\left\lceil #1\right\rceil}
\newcommand{\norm}[1]{\left\lVert #1\right\rVert}
\newcommand{\ip}[2]{\left\langle #1,#2\right\rangle}
\DeclareMathOperator{\diag}{diag}

\title{Fixed-space indices and symmetry breaking\\in variational problems}
\author{Daniele Garrisi\thanks{School of Mathematical Sciences, University of Nottingham Ningbo China, Ningbo, P. R. China.}
\and Alessandro Portaluri\thanks{DISAFA, Universit\`a degli Studi di Torino, Largo Paolo Braccini 2, 10095 Grugliasco, Torino, Italy.}
\and Li Wu\thanks{School of Mathematics, Shandong University, Jinan 250100, P. R. China.}}
\date{\today}
\begin{document}
\maketitle

\begin{abstract}
We establish a quantitative criterion for variational bifurcation
with prescribed exact stabilizers under finite group actions.
Combining ordinary spectral flow on fixed-point spaces with
minimal-overgroup nullities, we obtain a lower bound on the number
of exact-stabilizer bifurcation parameters when the restricted
singular set is finite. The theorem allows strongly indefinite
Fredholm Hessians and requires endpoint invertibility only on
the chosen fixed space, without parameter differentiability
or regular crossings. When the restricted Hessians have finite
Morse index, the criterion involves only indices and nullities.
We also treat varying form domains under explicit transport
assumptions and regular matrix Sturm--Liouville systems.

As our principal application, we consider a Newtonian nine-body
family consisting of two interlaced squares with fixed satellite
masses and a varying central mass. We prove analytically that
central configurations with two nonconjugate exact reflection
stabilizers bifurcate at a common mass in $(142,150)$, displacing
the central body from the centre of mass. The family contains
a reported counterexample to Moeckel's linear-stability conjecture;
our bifurcation proof is independent of that stability assertion.
\end{abstract}



\medskip
{\small\tableofcontents}

\section{Introduction}\label{sec:intro}

In the study of the Newtonian $N$-body problem, symmetry-breaking bifurcation of central configurations provides a natural mechanism for finding new families of relative equilibria. As the masses of the bodies vary, a symmetric central configuration may become degenerate; however, degeneracy of the Hessian alone neither forces bifurcation nor identifies the symmetries of the bifurcating solutions. Even when a change in the Morse index on the fixed-point subspace of a subgroup $H\le G$ yields non-trivial critical points by spectral-flow bifurcation on the fixed space \cite{FPR,PW} and Palais's principle of symmetric criticality \cite{Palais}, such critical points are only known a priori to have stabilizer containing $H$, and might in fact be fixed by a strictly larger subgroup $K>H$ (cf.\ \cite{SmollerWasserman,IJSW}). Therefore, if one wishes to show that the bifurcating solutions have {\em exact} stabilizer $H$, one has to take into account the behavior of the Hessian on the fixed-point subspaces of the overgroups of $H$.

In the present paper we address this issue from both an abstract and a mechanical viewpoint. First, in Theorem A we give a variational criterion which uses the spectral flow on $\cH^H$ together with the nullities of the Hessian on the minimal strict overgroups of $H$ in order to bound from below the number of bifurcation instants for critical points having exact stabilizer $H$. Second, in Theorem B we apply this criterion to a nine-body family formed by two concentric interlaced squares of satellites with fixed masses and a central body of varying mass $\mu(a)$. We show that central configurations having two non-conjugate exact reflection stabilizers in $D_4$ bifurcate at a {\em common critical mass} $\mu_*\in(142,150)$, with the central body displaced away from the center of mass along the corresponding axis of reflection. The fact that the two non-conjugate reflection classes bifurcate at the same critical mass is a direct consequence of the representation theory of $D_4$. More precisely, once the $16$-dimensional satellite block is shown to be positive definite on the considered interval, completing the square reduces the problem to a $2\times 2$ Schur complement on the central displacement plane $\R^2$; since this plane carries the irreducible standard representation of $D_4$, the Schur complement must be of the form $\sigma(a)I_2$, and hence the same scalar function $\sigma(a)$ governs both reflection-fixed subspaces.

It is worth noticing that the nine-body family considered here contains at $\mu=150$ the central configuration recently found numerically by Seamons \cite{Seamons2026} as a putative counterexample to Moeckel's conjecture (see \cite[Problem 16]{AlbouyCabralSantos}). By means of explicit rational interval estimates, we prove that as the central mass increases from $\mu_-\in(142,143)$ to $\mu=150$, the Morse index on the $14$-dimensional shape slice changes from $0$ to $2$ (and from $0$ to $1$ on each reflection-fixed slice), thus giving an analytic proof of the fact that the $\mu=150$ configuration is a saddle point of Morse index two on shape space, as well as of the symmetry-breaking bifurcation taking place along the branch.

\subsection{Fixed-space bifurcation and exact symmetry}
\label{subsec:intro_thmA}

Let us introduce the abstract variational setting. Let $I=[a,b]$, with $a<b$, and let $G$ be a finite group acting orthogonally on a real separable Hilbert space $\cH$. We consider a family of functionals
$$
 f_\lambda:B_\rho(0)\subset\cH\longrightarrow\R,
 \qquad \lambda\in I,
$$
and assume that the following conditions hold:
\begin{enumerate}[(H1)]
\item
For each $\lambda\in I$, $f_\lambda$ is of class $C^2$ on $B_\rho(0)$, and the mappings
$$
 (\lambda,u)\longmapsto f_\lambda(u),\qquad
 (\lambda,u)\longmapsto Df_\lambda(u)\in\cH^*,\qquad
 (\lambda,u)\longmapsto D^2f_\lambda(u)\in\mathcal L(\cH,\cH^*)
$$
are continuous on $I\times B_\rho(0)$ with respect to the norm topologies on $\cH^*$ and $\mathcal L(\cH,\cH^*)$.
\item
Each $f_\lambda$ is $G$-invariant on $B_\rho(0)$ and satisfies $Df_\lambda(0)=0$, $\forall\,\lambda\in I$.
\item
The bounded self-adjoint operator $L_\lambda\in\mathcal L(\cH)$ defined by
$$
 D^2f_\lambda(0)[v,w]=\ip{L_\lambda v}{w}, \quad \forall\, v,w\in\cH,
$$
is a Fredholm operator for every $\lambda\in I$.
\end{enumerate}

Given any subgroup $K\le G$, we set
$$
 \cH^K=\big\{u\in\cH\ \big|\ ku=u,\ \forall\,k\in K\big\},
 \qquad
 L_\lambda^K=L_\lambda\big|_{\cH^K}\in\mathcal L(\cH^K),
 \qquad
 \nu_K(t)=\dim\ker L_t^K.
$$
Since $L_\lambda$ is $G$-equivariant, $\cH^K$ is an $L_\lambda$-invariant reducing subspace and $\lambda\mapsto L_\lambda^K$ defines a norm-continuous path of bounded self-adjoint Fredholm operators on $\cH^K$. Furthermore, by Palais's Principle of Symmetric Criticality \cite{Palais}, $u\in B_\rho(0)\cap\cH^K$ is a critical point of $f_\lambda|_{\cH^K}$ if and only if it is a critical point of $f_\lambda$ on $\cH$.

Let now $H\le G$ be a fixed subgroup such that $L_a^H$ and $L_b^H$ be invertible on $\cH^H$ (whereas the ambient operators $L_a$ and $L_b$ on $\cH$ may well be degenerate). We denote by
\[
 s_H=\sflo(L^H;[a,b])\in\mathbb Z
\]
the spectral flow of the restricted path $L^H=\{L_\lambda^H\}_{\lambda\in[a,b]}$, where eigenvalues crossing zero from $(-\infty,0)$ to $[0,\infty)$ are counted with a positive sign \cite{Phillips1996,BLP,PW}. In the special case in which $L_\lambda^H$ has finite Morse index $m^-(L_\lambda^H)<\infty$ for all $\lambda\in I$, we simply have
\begin{equation}\label{eq:intro_morse}
 s_H=m^-(L_a^H)-m^-(L_b^H).
\end{equation}

To separate critical points with exact stabilizer $H$ from those having a larger symmetry group, we introduce the family of minimal strict overgroups of $H$ in $G$, namely
\[
 \mathcal K_H=
 \big\{K\le G\ \big|\ H<K,\text{ and there is no }J\le G\text{ with }H<J<K\big\}.
\]
Because $G$ is a finite group, any subgroup $S\le G$ strictly containing $H$ contains at least one subgroup $K\in\mathcal K_H$, which implies $\cH^S\subset\cH^K$. In other words, in order to rule out stabilizers strictly larger than $H$ near a given instant $t\in(a,b)$, it is enough to check that no non-trivial critical points bifurcate in the spaces $\cH^K$ for $K\in\mathcal K_H$. Accordingly, we set
\begin{equation}\label{eq:exceptional}
 \begin{aligned}
 \Sigma_H&=\big\{t\in(a,b)\ \big|\ \nu_H(t)>0\big\},\\
 \cE_H&=\big\{t\in\Sigma_H\ \big|\ \nu_K(t)>0\text{ for some }K\in\mathcal K_H\big\}.
 \end{aligned}
\end{equation}
If $t\in\Sigma_H\setminus\cE_H$, then $\nu_K(t)=\dim\ker L_t^K=0$ for every $K\in\mathcal K_H$; thus each $L_t^K$ is invertible on $\cH^K$, and the trivial solution $u=0$ is locally isolated in $\cH^K$ for $\lambda$ close to $t$. Whenever $\Sigma_H$ is a finite set, we define the integers
\begin{equation}\label{eq:BH_MH}
 B_H=\sum_{t\in\cE_H}\nu_H(t)
 \quad \textrm{ and } \quad
 M_H=\max\big(\{\nu_H(t):t\in\Sigma_H\setminus\cE_H\}\cup\{0\}\big),
\end{equation}
where as usual an empty sum is understood to be zero.

\begin{theoremA}[Bifurcation and prescribed exact stabilizers]
Suppose that \textup{(H1)--(H3)} hold, and let $H\le G$ be a subgroup such that $L_a^H$ and $L_b^H$ be invertible on $\cH^H$. Then the following hold:
\begin{enumerate}[(i)]
\item
If $s_H\ne0$, there exist $t\in(a,b)$ and sequences $\lambda_j\to t$, $u_j\in\cH^H\setminus\{0\}$ with $u_j\to0$, such that
$$
 \nabla f_{\lambda_j}(u_j)=0, \quad \forall\, j\ge1.
$$
\item
If $\Sigma_H$ is finite and $|s_H|>B_H$, then $M_H\ge1$ and the number $N_H$ of distinct instants $t\in\Sigma_H\setminus\cE_H$ admitting sequences $(\lambda_j,u_j)\to(t,0)$ of critical points with $\Stab_G(u_j)=H$ satisfies
\begin{equation}\label{eq:thmA_bound}
 N_H\ge\ceil{\frac{|s_H|-B_H}{M_H}}.
\end{equation}
\end{enumerate}
In particular, whenever $m^-(L_\lambda^H)<\infty$ on $I$, $s_H$ can be computed via \eqref{eq:intro_morse}.
\end{theoremA}

The idea behind the proof of Theorem A can be summarized as follows. At each isolated singular instant $t\in\Sigma_H$, the local spectral flow $s_H(t)=\sflo(L^H;[t-\delta,t+\delta])$ satisfies $|s_H(t)|\le\nu_H(t)$ by Lemma~\ref{lem:localflow}. By path additivity of the spectral flow, the total contribution coming from the exceptional set $\cE_H$ is at most $B_H$. Therefore, under the assumption that $|s_H|>B_H$, the remaining spectral flow $|s_H|-B_H>0$ must be accounted for by non-exceptional instants $t\in\Sigma_H\setminus\cE_H$ for which $s_H(t)\ne0$. At each such instant $t$, the classical variational bifurcation theorem \cite{FPR,PW} applied on $\cH^H$ gives a sequence of non-trivial $H$-invariant critical points bifurcating from $(t,0)$, whereas the invertibility of $L_t^K$ on $\cH^K$ ($K\in\mathcal K_H$) rules out any stabilizer strictly larger than $H$.

\subsection{Simultaneous bifurcation in the nine-body problem}
\label{subsec:intro_thmB}

We now pass to describe our main result for the planar nine-body problem. The configurations we consider are formed by four particles of mass $1$ located at $(\pm1,0)$ and $(0,\pm1)$, four particles of mass $\beta_0>0$ located at $(\pm a,\pm a)$, and a central particle of mass $\mu(a)$ placed at the origin; see \eqref{eq:ninepositions}--\eqref{eq:ninemasses} for the explicit expressions of $\beta_0$ and $\mu(a)$. While the eight satellite masses are kept fixed, the geometric parameter $a$ varies in the closed rational interval
$$
 \mathcal J=[a_0,a_1]
 =\left[\frac{177}{250},\frac{14161}{20000}\right],
$$
on which the central mass $\mu(a)$ is strictly decreasing from $\mu(a_0)=150$ to $\mu(a_1)=\mu_-\in(142,143)$.

For each $a\in\mathcal J$, this concentric central configuration is invariant under the action of the dihedral group $D_4=\langle r,s\rangle$ of order eight, where $r$ acts on the plane as a counterclockwise rotation by $\pi/2$ and $s$ as the reflection across the horizontal axis, accompanied by the corresponding permutations of equal-mass satellites. Up to conjugacy in $D_4$, we have two distinct reflection subgroups of order two, namely the axial reflection and the diagonal reflection:
$$
 H_{\mathrm{ax}}=\langle s\rangle
 \quad \textrm{ and } \quad
 H_{\mathrm{diag}}=\langle rs\rangle.
$$
Once translations, dilations, and planar rotations are quotiented out by passing to a local $D_4$-equivariant $14$-dimensional shape slice $W$, we obtain the following theorem. Notice that, in contrast with bifurcations taking place within concentric rosette configurations \cite{LeiSantoprete}, the bifurcating solutions found here move the central mass away from the origin.

\begin{theoremB}[Simultaneous reflection bifurcation]
For the Newtonian nine-body family \eqref{eq:ninepositions}--\eqref{eq:ninemasses}, the following statements hold:
\begin{enumerate}[(i)]
\item
As $\mu$ increases from $\mu_-$ to $150$, the Morse index on each $7$-dimensional reflection-fixed shape subspace $W^{H_{\mathrm{ax}}}$ resp.\ $W^{H_{\mathrm{diag}}}$ changes from $0$ to $1$, whereas the Morse index on the full $14$-dimensional shape space $W$ changes from $0$ to $2$. Accordingly, the corresponding spectral flows are $-1$, $-1$, and $-2$.
\item
There exists a common critical mass $\mu_*\in(\mu_-,150)\subset(142,150)$ at which central configurations with exact stabilizer $H_{\mathrm{ax}}$ and central configurations with exact stabilizer $H_{\mathrm{diag}}$ bifurcate simultaneously from the concentric branch, modulo planar rotations and scaling.
\item
Along both bifurcating sequences, any sufficiently close non-trivial configuration displaces the central body ($q_0\ne0$) away from the center of mass along the respective axis of reflection.
\end{enumerate}
\end{theoremB}

By saying that the two bifurcations occur simultaneously at $\mu_*$, we mean that both reflection classes admit sequences of non-trivial central configurations accumulating at the concentric branch at the same limiting mass $\mu_*$ (of course, the discrete mass values along the two sequences need not be equal term by term; establishing that $\mu_*$ is unique and that the bifurcating sets form smooth curves requires higher-order derivative estimates, which we discuss in Section 7).

The reason why both reflection types bifurcate at the very same critical mass $\mu_*$ can be easily seen from the block structure of the rotation-lifted Hessian $\mathscr Q(a)$ (defined in Lemma~\ref{lem:mechanicalinertia}). If we split the configuration space $\R^{18}=\R^2\oplus\R^{16}$ into the central-body displacement $\xi_0\in\R^2$ and the satellite displacements $\xi_{\mathrm{sat}}\in\R^{16}$, then $\mathscr Q(a)$ reads
\[
 \mathscr Q(a)=\begin{pmatrix}\alpha(a)I_2 & B(a)^T \\ B(a) & C(a)\end{pmatrix}.
\]
In Lemma~\ref{lem:ninecertificate} we decompose $\R^{16}$ into $D_4$-isotypic blocks and use outward-rounded rational interval arithmetic to show that the $16\times16$ satellite block $C(a)=\mathscr Q(a)|_{\{\xi_0=0\}}$ is strictly positive definite for all $a\in\mathcal J$. From the positive definiteness of $C(a)$ we deduce two facts:
\begin{itemize}
\item[-] Every subgroup $K\le D_4$ strictly containing $H_{\mathrm{ax}}$ or $H_{\mathrm{diag}}$ contains the central inversion $r^2$, which acts as $-I_2$ on $\xi_0\in\R^2$. Hence $(\R^{18})^K\subset\{\xi_0=0\}$, so that $\mathscr Q(a)|_{(\R^{18})^K}$ is positive definite; in particular, $\nu_K(a)=0$ for all $a\in\mathcal J$ (which gives $B_{H_{\mathrm{ax}}}=B_{H_{\mathrm{diag}}}=0$), and no non-trivial central configuration with $q_0=0$ can bifurcate from the branch.
\item[-] Completing the square with respect to $\xi_{\mathrm{sat}}$ reduces the quadratic form to the $2\times2$ Schur complement
$$
 S(a)=\alpha(a)I_2-B(a)^TC(a)^{-1}B(a)
$$
on $\R^2$. Since $D_4$ acts irreducibly on $\R^2$ via the standard representation $E$, Schur's lemma implies that $S(a)=\sigma(a)I_2$ for a real-analytic scalar function $\sigma(a)$.
\end{itemize}
As a consequence, the restrictions of $S(a)$ to the axial line $(\R^2)^{H_{\mathrm{ax}}}=\R(1,0)$ and to the diagonal line $(\R^2)^{H_{\mathrm{diag}}}=\R(1,1)$ are both given by multiplication by the same scalar $\sigma(a)$. Since our endpoint estimates show that $\sigma(a_0)<0<\sigma(a_1)$, any zero $a_*\in(a_0,a_1)$ of $\sigma(a)$ at which the local spectral flow is non-zero produces bifurcation in $W^{H_{\mathrm{ax}}}$ and $W^{H_{\mathrm{diag}}}$ at the same time.

\paragraph{Connection with Moeckel's conjecture.}
Moeckel's conjecture \cite[Problem 16]{AlbouyCabralSantos} asks whether a linearly stable relative equilibrium in the planar $N$-body problem must come from a central configuration which is a non-degenerate minimum of the potential on the shape sphere (i.e.\ with shape Morse index zero). In a uniformly rotating frame with angular speed $\omega(a)$, the linearized Hamiltonian system on the $14$-dimensional shape slice $W$ is associated with the quadratic matrix pencil
$$
 P_W(z)=z^2M_W+2\omega(a)z(MJ)_W-K_W(a),
$$
where $K_W(a)=(D^2U(q(a))+\omega(a)^2M(a))|_W$ denotes the Hessian of the amended potential on $W$. By a classical result of Hu and Sun \cite{HuSun2009}, if a non-degenerate central configuration has odd Morse index on $W$, then $\det P_W(0)=\det(-K_W)<0$ (as $\dim W=2N-4$ is even), and hence $\det P_W(z)=0$ has at least one positive real root, which implies spectral instability. When $m^-(K_W)$ is even and positive, instead, $\det P_W(0)>0$, so that gyroscopic stabilization due to the Coriolis term $2\omega(a)z(MJ)_W$ is not ruled out by the determinant parity. Theorem B shows analytically that for the Seamons configuration at $\mu(a_0)=150$ \cite{Seamons2026} one has $m^-(K_W(a_0))=2$, whereas $m^-(K_W(a_1))=0$ at $\mu_-\in(142,143)$. At least one intermediate mass supports simultaneous bifurcation of the two exact reflection types.

\subsection{Relation to existing equivariant bifurcation theories}
\label{subsec:intro_literature}

Let us briefly compare Theorem A with related results in equivariant bifurcation theory. Part (i) of Theorem A follows by applying the spectral-flow bifurcation theorem of Fitzpatrick, Pejsachowicz, and Recht \cite{FPR} (in the version for continuous $C^2$ families proved by Pejsachowicz and Waterstraat \cite{PW}) to $f_\lambda|_{\cH^H}$, together with Palais's principle \cite{Palais}. Part (ii) combines this existence result with a local isolation argument on the overgroup fixed spaces $\cH^K$, $K\in\mathcal K_H$.

In the recent papers \cite{IJW2021,IJSW,IJSW2026}, Izydorek, Janczewska, Starostka, and Waterstraat constructed an $RO(G)$-valued equivariant spectral flow $\sflo_G(L;[a,b])$ for paths of $G$-equivariant self-adjoint Fredholm operators having invertible endpoints $L_a,L_b$ on $\cH$. Whenever $L_a$ and $L_b$ are invertible on $\cH$, the fixed-point dimension homomorphism $\Phi_H([U]-[V])=\dim U^H-\dim V^H$ maps $\sflo_G(L;[a,b])$ onto $s_H$. At the linear level, $\sflo_G(L;[a,b])\in RO(G)$ contains strictly more information than the family of integer flows $\{s_K\}_{K\le G}$: for example, if $G=C_5$, the two non-isomorphic $2$-dimensional irreducible rotation representations $E_1$ and $E_2$ satisfy $\dim E_1^K=\dim E_2^K$ for all $K\le C_5$, so that the path $L_t=tI_{E_1}\oplus(-tI_{E_2})$ on $[-1,1]$ has $\sflo_{C_5}(L;[-1,1])=[E_1]-[E_2]\ne0$ in $RO(C_5)$, whereas $s_K=0$ for all $K\le C_5$ (cf.\ Remark~\ref{rem:comparison}). On the other hand, working with the integer flow $s_H$ and the overgroup nullities $\nu_K(t)$ only requires endpoint invertibility on the subspace $\cH^H$ (allowing $L_a$ and $L_b$ to have non-trivial kernels on $\cH$), and directly yields the lower bound \eqref{eq:thmA_bound} on the number of bifurcation instants with exact stabilizer $H$.

A different approach is provided by equivariant degree theories, such as the orthogonal degree of Ize and Vignoli \cite{IzeVignoli1999,IzeVignoliBook} (used by Garc\'ia-Azpeitia and Ize \cite{GarciaAzpeitiaIze2011,GarciaAzpeitiaIze} in the study of $N$-body relative equilibria and periodic orbits), Kowalczyk's degree for invariant potentials \cite{Kowalczyk}, and twisted equivariant degrees \cite{BKRSteinlein}. While equivariant degree methods provide global connected branches for compact perturbations of the identity via Burnside-ring algebra, Theorem A is a local variational criterion applicable to general norm-continuous Fredholm paths using only restricted Morse indices and nullities.

\subsection{Organization of the paper}
\label{subsec:intro_outline}

We finish this introduction with a short summary of the content of the paper. In Section 2 we recall the basic facts on fixed-point subspaces and spectral flow, prove the budget estimate (Theorem~\ref{thm:budget}) and Theorem A, review the equivariant Lyapunov--Schmidt reduction and simple-crossing corollary, and discuss the comparison with $RO(G)$-valued spectral flow. In Section 3 we deal with families of unbounded self-adjoint operators $A_\lambda$ having varying form domains $V_\lambda=D(|A_\lambda|^{1/2})$, showing how the canonical isometries $U_\lambda=(1+A_\lambda^2)^{-1/4}:\cH\to V_\lambda$ allow one to pull back the quadratic forms to the fixed space $\cH$ (Theorem~\ref{thm:transport}). Section 4 is devoted to regular matrix Sturm--Liouville problems and moving boundary conditions, including the Robin-to-Dirichlet example where Riesz continuity fails. In Section 5 we give the version of Theorem A for equivariant Lagrangian intersections in a Weinstein chart. In Section 6 we turn to the Newtonian $N$-body problem: we prove the rotational lift lemma (Lemma~\ref{lem:mechanicalinertia}), revisit Palmore's centered equilateral triangle bifurcation, carry out the $D_4$-block decomposition and rational interval estimates needed to prove Theorem B for the nine-body family, and study spatial relative-periodic orbits of a vertically oscillating ring. Finally, in Section 7 we discuss two open questions on the nine-body family concerning the fourth-order Lyapunov--Schmidt normal form and the degree-$14$ gyroscopic polynomial.

\section{Fixed-space bifurcation and exact stabilizers}\label{sec:symmetry}

The goal of this section is to prove Theorem A on a fixed real separable Hilbert space $\cH$. Throughout the section we take $I=[a,b]$, with $a<b$, and assume that $G$ is a finite group acting orthogonally on $\cH$ and that $f_\lambda:B_\rho(0)\subset\cH\to\R$ satisfies hypotheses (H1)--(H3) of Section~\ref{subsec:intro_thmA}.

\subsection{Fixed-point subspaces and the classical bifurcation input}\label{subsec:classical}

Passing to the fixed-point subspace $\cH^H$ of a subgroup $H\le G$ has a double advantage: on the one hand, critical points in $\cH^H$ automatically possess symmetry at least $H$; on the other hand, restricting to $\cH^H$ can kill undesired kernels present in the full space $\cH$ (such as zero modes produced by continuous symmetries or other symmetry modes at the endpoints). Let us briefly check that self-adjointness, the Fredholm property, norm continuity, and criticality are all preserved when restricting to $\cH^K$ for $K\le G$.

For any subgroup $K\le G$, we define the bounded linear operator
$$
 P_K=\frac{1}{|K|}\sum_{k\in K}k\in\mathcal{L}(\cH).
$$
Since the action of $G$ is orthogonal ($k^*=k^{-1}$) and $K$ is a subgroup, it is readily seen that $P_K^*=P_K$, $P_K^2=P_K$, and $\ran P_K=\cH^K$; thus $P_K$ is the orthogonal projection onto $\cH^K$. Moreover, differentiating the identity $f_\lambda(gu)=f_\lambda(u)$ twice at $u=0$, we see that $L_\lambda$ commutes with the action of $G$, namely
$$
 L_\lambda g=gL_\lambda, \quad \forall\, g\in G.
$$
It follows that $P_KL_\lambda=L_\lambda P_K$, so that both $\cH^K$ and $(\cH^K)^\perp$ are $L_\lambda$-invariant reducing subspaces. Therefore, the restriction $L_\lambda^K=L_\lambda|_{\cH^K}\in\mathcal{L}(\cH^K)$ is a bounded self-adjoint operator satisfying
$$
 \ker L_\lambda^K=\ker L_\lambda\cap\cH^K=(\ker L_\lambda)^K
 \quad \textrm{ and } \quad
 \nu_K(\lambda):=\dim\ker L_\lambda^K<\infty.
$$
Also, since $L_\lambda$ is Fredholm, its range $\ran L_\lambda$ is closed in $\cH$, which implies that $P_K(\ran L_\lambda)=\ran L_\lambda^K=\ran L_\lambda\cap\cH^K$ is closed in $\cH^K$. This shows that $\lambda\mapsto L_\lambda^K$ is a norm-continuous path of bounded self-adjoint Fredholm operators on $\cH^K$.

At the non-linear level, we have Palais's Principle of Symmetric Criticality \cite{Palais}, whose simple proof is included for the reader's convenience.

\begin{theorem}[Principle of Symmetric Criticality \cite{Palais}]\label{thm:PSC}
Let a compact group $G$ act orthogonally on a real Hilbert space $\cH$, and let $f:U\to\R$ be a $G$-invariant $C^1$ functional on a $G$-invariant open set $U\subset\cH$. For any closed subgroup $K\le G$ and any $u\in U\cap\cH^K$, $u$ is a critical point of the restriction $f|_{U\cap\cH^K}$ if and only if $\nabla f(u)=0$ in $\cH$.
\end{theorem}
\begin{proof}
Orthogonal invariance $f(ku)=f(u)$ differentiates to $\nabla f(ku)=k\nabla f(u)$ for all $k\in K$. If $u\in\cH^K$, then $ku=u$, so $k\nabla f(u)=\nabla f(u)$ for all $k\in K$. Therefore, $\nabla f(u)\in\cH^K$. If $u$ is a critical point of $f$ constrained to $\cH^K$, then $\nabla f(u)\perp\cH^K$. Since $\nabla f(u)\in\cH^K\cap(\cH^K)^\perp$, $u$ is a non-constrained critical point.
\end{proof}

We adopt the definition of spectral flow given by Phillips \cite{Phillips1996}, Booss-Bavnbek, Lesch, and Phillips \cite{BLP}, and Pejsachowicz and Waterstraat \cite{PW}. For a norm-continuous path $L=\{L_\lambda\}_{\lambda\in[a,b]}$ of bounded self-adjoint Fredholm operators on a real separable Hilbert space, one chooses a partition $a=\lambda_0<\lambda_1<\dots<\lambda_N=b$ and numbers $\epsilon_i>0$ such that $\pm\epsilon_i\notin\spec(L_\lambda)$ and $[-\epsilon_i,\epsilon_i]\cap\spec_{\mathrm{ess}}(L_\lambda)=\varnothing$ for every $\lambda\in[\lambda_{i-1},\lambda_i]$. Then the spectral projections $\mathbf{1}_{[-\epsilon_i,\epsilon_i]}(L_\lambda)$ depend norm-continuously on $\lambda\in[\lambda_{i-1},\lambda_i]$ and have finite rank, and the spectral flow of $L$ on $[a,b]$ is given by
\[
 \sflo(L;[a,b])=\sum_{i=1}^N\Big(\dim\ran\mathbf{1}_{[0,\epsilon_i]}(L_{\lambda_i})-\dim\ran\mathbf{1}_{[0,\epsilon_i]}(L_{\lambda_{i-1}})\Big)\in\mathbb{Z}.
\]
With this sign convention, eigenvalues crossing zero from $(-\infty,0)$ to $[0,\infty)$ give a positive contribution. $\sflo(L;[a,b])$ does not depend on the chosen partition or on the numbers $\epsilon_i$, is additive with respect to concatenation of paths and direct sums, vanishes if $L_\lambda$ is invertible for every $\lambda\in[a,b]$, and is homotopy invariant within the class of paths having invertible endpoints. In particular, when $m^-(L_\lambda)=\dim\ran\mathbf{1}_{(-\infty,0)}(L_\lambda)<\infty$ for all $\lambda\in[a,b]$ and the endpoints $L_a,L_b$ are invertible, one can use a single spectral window for the whole path (or telescope the sum above) to get
\begin{equation}\label{eq:sfmorse}
 \sflo(L;[a,b])=m^-(L_a)-m^-(L_b).
\end{equation}
When $L_a,L_b$ are invertible and $t\in(a,b)$ is an isolated singular instant of $L$ (meaning that $L_\lambda$ is invertible for $0<|\lambda-t|\le\delta$), the local spectral flow
$$
 \sflo(L;t):=\sflo(L;[t-\delta,t+\delta])
$$
is well-defined and independent of $\delta>0$ small enough. We have the following simple estimate on $\sflo(L;t)$.

\begin{lemma}[Local spectral-flow bound]\label{lem:localflow}
Let $L=\{L_\lambda\}_{\lambda\in[a,b]}$ be a norm-continuous path of bounded self-adjoint Fredholm operators, and let $t\in(a,b)$ be an isolated singular parameter of $L$. Then
\begin{equation}\label{eq:localflow_bound}
 |\sflo(L;t)|\le\dim\ker L_t.
\end{equation}
\end{lemma}
\begin{proof}
Since $L_t$ is self-adjoint and Fredholm, $0$ is an isolated point of $\spec(L_t)$ and an eigenvalue of finite multiplicity $d=\dim\ker L_t$. Choose $\epsilon>0$ such that $\spec(L_t)\cap[-\epsilon,\epsilon]=\{0\}$. By norm continuity of $\lambda\mapsto L_\lambda$, there exists $\delta>0$ such that $[t-\delta,t+\delta]\subset(a,b)$, $L_\lambda$ is invertible for all $\lambda\in[t-\delta,t+\delta]\setminus\{t\}$, and $\pm\epsilon\notin\spec(L_\lambda)$ with $[-\epsilon,\epsilon]\cap\spec_{\mathrm{ess}}(L_\lambda)=\varnothing$ for all $\lambda\in[t-\delta,t+\delta]$. The finite-rank spectral projection $P(\lambda)=\mathbf{1}_{[-\epsilon,\epsilon]}(L_\lambda)$ is then norm-continuous on $[t-\delta,t+\delta]$ and therefore has constant rank $\dim\ran P(\lambda)=\dim\ran P(t)=d$. At the invertible endpoints $t_\pm=t\pm\delta$, $0$ is not an eigenvalue of $L_{t_\pm}$, so $k_\pm:=\dim\ran\mathbf{1}_{[0,\epsilon]}(L_{t_\pm})$ satisfies $0\le k_\pm\le d$. By definition, $\sflo(L;t)=k_+-k_-$, whence $|\sflo(L;t)|=|k_+-k_-|\le d=\dim\ker L_t$.
\end{proof}

Let us now recall the variational bifurcation theorem of Fitzpatrick, Pejsachowicz, and Recht \cite{FPR} in the formulation for continuous $C^2$ families due to Pejsachowicz and Waterstraat \cite[Theorem 2.1]{PW} (cf.\ also \cite{Pejsachowicz2025}).

\begin{theorem}[Classical variational spectral-flow bifurcation {\cite[Theorem 2.1]{PW}}]\label{thm:classical}
Let $f_\lambda:B_\rho(0)\subset\cH\to\R$, $\lambda\in[a,b]$, be a family of functionals satisfying \textup{(H1)--(H3)} (with trivial group), so that $\nabla f_\lambda(0)=0$ and $L_\lambda=D^2f_\lambda(0)$ is a norm-continuous path of bounded self-adjoint Fredholm operators.
\begin{enumerate}[(i)]
\item If $L_a$ and $L_b$ are invertible and $\sflo(L;[a,b])\ne0$, then there exist $t\in(a,b)$, a sequence $\lambda_j\to t$, and $u_j\in B_\rho(0)\setminus\{0\}$ such that $u_j\to0$ and $\nabla f_{\lambda_j}(u_j)=0$.
\item In particular, if $t\in(a,b)$ is an isolated singular instant of $L$ with $\sflo(L;t)\ne0$, then $t$ is a bifurcation instant of critical points.
\item If $L_a,L_b$ are invertible, the singular set $\Sigma=\{t\in(a,b):\dim\ker L_t>0\}$ is finite, and $\sflo(L;[a,b])\ne0$, then $\Sigma\ne\varnothing$ and the number of distinct bifurcation instants in $(a,b)$ is at least
$$
 \ceil{\frac{|\sflo(L;[a,b])|}{\max_{t\in\Sigma}\dim\ker L_t}}.
$$
\end{enumerate}
\end{theorem}

We will also need the following standard local isolation property near any parameter at which the Hessian is invertible.

\begin{lemma}[Uniform local isolation at an invertible Hessian]\label{lem:uniqueness}
Under hypotheses \textup{(H1)--(H3)} on a Hilbert space $\cW$, suppose that $L_{t_0}=D^2f_{t_0}(0)\in\mathcal{L}(\cW)$ be invertible at some $t_0\in[a,b]$. Then there exist $\delta>0$ and $r\in(0,\rho)$ such that whenever $\lambda\in[a,b]$ satisfies $|\lambda-t_0|<\delta$ and $u\in\cW$ satisfies $\norm{u}<r$ and $\nabla f_\lambda(u)=0$, one has $u=0$.
\end{lemma}
\begin{proof}
Set $R=L_{t_0}^{-1}\in\mathcal{L}(\cW)$. By the joint norm continuity of $(\lambda,u)\mapsto D^2f_\lambda(u)$ at $(t_0,0)$, we can choose $\delta>0$ and $r\in(0,\rho)$ such that $\norm{I-RD^2f_\lambda(u)}=\norm{R(L_{t_0}-D^2f_\lambda(u))}\le\frac{1}{2}$ for all $\lambda\in[a,b]$ with $|\lambda-t_0|<\delta$ and all $u\in B_r(0)\subset\cW$. For each such $\lambda$, the map $\Phi_\lambda(u)=u-R\nabla f_\lambda(u)$ satisfies $\Phi_\lambda(0)=0$ and $\norm{D\Phi_\lambda(u)}\le\frac{1}{2}$ on $B_r(0)$, which implies $\norm{\Phi_\lambda(u)}\le\frac{1}{2}\norm{u}$. Since any critical point $u\in B_r(0)$ of $f_\lambda$ is a fixed point of $\Phi_\lambda$, we obtain $\norm{u}\le\frac{1}{2}\norm{u}$, and hence $u=0$.
\end{proof}

\subsection{The exclusion estimate and proof of Theorem A}\label{subsec:exact}

We are now in a position to put together Theorem~\ref{thm:classical}, Theorem~\ref{thm:PSC}, Lemma~\ref{lem:localflow}, and Lemma~\ref{lem:uniqueness} to prove the exact-stabilizer budget theorem, from which Theorem A will follow at once.

Recall from \eqref{eq:exceptional}--\eqref{eq:BH_MH} that for a fixed subgroup $H\le G$, $\mathcal{K}_H$ denotes the set of minimal strict overgroups of $H$ in $G$, $\Sigma_H=\{t\in(a,b):\nu_H(t)>0\}$, and $\cE_H=\{t\in\Sigma_H:\nu_K(t)>0\text{ for some }K\in\mathcal{K}_H\}$. Notice that each parameter $t\in\cE_H$ is counted once in the sum $B_H=\sum_{t\in\cE_H}\nu_H(t)$, regardless of how many overgroups $K\in\mathcal{K}_H$ have a non-trivial kernel at $t$. In the extreme case $H=G$, one has $\mathcal{K}_G=\varnothing$, so that $\cE_G=\varnothing$ and $B_G=0$. When $\Sigma_H$ is finite and $L_a^H,L_b^H$ are invertible, every $t\in\Sigma_H$ is an isolated singular parameter of $L^H$, and we write $s_H(t)=\sflo(L^H;t)$ for the local spectral flow of $L^H$ across $t$.

\begin{theorem}[Exact-stabilizer budget and parameter lower bound]\label{thm:budget}
Assume \textup{(H1)--(H3)} hold, fix $H\le G$, and suppose that $L_a^H,L_b^H$ be invertible on $\cH^H$ and $\Sigma_H$ be finite. Define the active non-exceptional crossing set
\begin{equation}\label{eq:TH_def}
 \mathcal{T}_H=\big\{t\in\Sigma_H\setminus\cE_H\ \big|\ s_H(t)\ne0\big\}.
\end{equation}
\begin{enumerate}[(i)]
\item For every $t\in\mathcal{T}_H$, there exist a sequence $\lambda_j\to t$ and non-zero critical points $u_j\in\cH^H\setminus\{0\}$ of $f_{\lambda_j}$ on $\cH$ such that $u_j\to0$ and $\Stab_G(u_j)=H$ for all $j\ge1$.
\item The nullities on $\mathcal{T}_H$ satisfy
\begin{equation}\label{eq:count}
 \sum_{t\in\mathcal{T}_H}\nu_H(t)\ge |s_H|-B_H.
\end{equation}
\item As a consequence, whenever $|s_H|>B_H$, we have $\mathcal{T}_H\ne\varnothing$, $M_H=\max_{t\in\Sigma_H\setminus\cE_H}\nu_H(t)\ge1$, and the number $N_H$ of distinct exact-$H$ bifurcation instants in $(a,b)$ satisfies
\begin{equation}\label{eq:NH_ceiling}
 N_H\ge\#\mathcal{T}_H\ge\ceil{\frac{|s_H|-B_H}{M_H}}.
\end{equation}
\end{enumerate}
\end{theorem}
\begin{proof}
\emph{Step 1: Localization of the spectral flow and proof of \eqref{eq:count}.}
Since $\Sigma_H\subset(a,b)$ is finite and $L_a^H,L_b^H$ are invertible, around each $t\in\Sigma_H$ we can pick pairwise disjoint closed intervals $J_t=[t-\eta_t,t+\eta_t]\subset(a,b)$ such that $J_t\cap\Sigma_H=\{t\}$. Outside the union $\bigcup_{t\in\Sigma_H}(t-\eta_t,t+\eta_t)$, the restricted path $L^H$ is invertible and hence its spectral flow vanishes. By path additivity of the spectral flow on $\cH^H$, we get
$$
 s_H=\sflo(L^H;[a,b])=\sum_{t\in\Sigma_H}\sflo(L^H;J_t)=\sum_{t\in\Sigma_H}s_H(t).
$$
By Lemma~\ref{lem:localflow} applied to $L^H$ at each $t\in\Sigma_H$, we have $|s_H(t)|\le\dim\ker L_t^H=\nu_H(t)$. Since $\Sigma_H$ is the disjoint union of $\cE_H$, $\mathcal{T}_H$, and $\{t\in\Sigma_H\setminus\cE_H:s_H(t)=0\}$, the triangle inequality gives
$$
 |s_H|\le\sum_{t\in\cE_H}|s_H(t)|+\sum_{t\in\mathcal{T}_H}|s_H(t)|\le\sum_{t\in\cE_H}\nu_H(t)+\sum_{t\in\mathcal{T}_H}\nu_H(t)=B_H+\sum_{t\in\mathcal{T}_H}\nu_H(t),
$$
proving \eqref{eq:count}.

\emph{Step 2: Existence of $H$-fixed critical points near each $t\in\mathcal{T}_H$.}
Let $t\in\mathcal{T}_H$ be fixed. Because $s_H(t)=\sflo(L^H;J_t)\ne0$ and $L^H$ is invertible at $\partial J_t$, we can apply Theorem~\ref{thm:classical}(i) to $f_\lambda|_{\cH^H}$ on $J_t$ to find $t'\in(t-\eta_t,t+\eta_t)$ and a sequence $(\lambda_j,u_j)\in J_t\times(\cH^H\setminus\{0\})$ converging to $(t',0)$ such that $\nabla(f_{\lambda_j}|_{\cH^H})(u_j)=0$. By Lemma~\ref{lem:uniqueness} on $\cW=\cH^H$, any bifurcation instant in $J_t$ must lie in $\Sigma_H\cap J_t=\{t\}$, which forces $t'=t$. By Theorem~\ref{thm:PSC}, each $u_j$ is a free critical point ($\nabla f_{\lambda_j}(u_j)=0$ in $\cH$), and $u_j\in\cH^H$ implies $\Stab_G(u_j)\ge H$.

\emph{Step 3: Exclusion of strictly larger stabilizers.}
For $t\in\mathcal{T}_H\subset\Sigma_H\setminus\cE_H$, we have $\nu_K(t)=\dim\ker L_t^K=0$ for all $K\in\mathcal{K}_H$. Since $L_t^K$ is self-adjoint and Fredholm with trivial kernel, $L_t^K\in\mathcal{L}(\cH^K)$ is invertible. By Lemma~\ref{lem:uniqueness} applied to $f_\lambda|_{\cH^K}$ on $\cW=\cH^K$ at $t_0=t$, there exist $\delta_K>0$ and $r_K\in(0,\rho)$ such that $u=0$ is the only critical point of $f_\lambda|_{\cH^K}$ satisfying $|\lambda-t|<\delta_K$ and $\norm{u}<r_K$.

If $\mathcal{K}_H=\varnothing$ (i.e.\ $H=G$), then $\Stab_G(u_j)\ge G$ directly gives $\Stab_G(u_j)=G=H$. Otherwise, since $G$ is finite, $\mathcal{K}_H$ is a finite non-empty set, and we can set $\delta_*:=\min_{K\in\mathcal{K}_H}\delta_K>0$ and $r_*:=\min_{K\in\mathcal{K}_H}r_K>0$. Since $(\lambda_j,u_j)\to(t,0)$, we have $|\lambda_j-t|<\delta_*$ and $\norm{u_j}<r_*$ for all $j\ge j_0$ large enough. If $S_j:=\Stab_G(u_j)$ strictly contained $H$ for some $j\ge j_0$, then $S_j$ would contain at least one minimal strict overgroup $K\in\mathcal{K}_H$, so that $u_j\in\cH^{S_j}\subset\cH^K$ would be a critical point of $f_{\lambda_j}|_{\cH^K}$ with $|\lambda_j-t|<\delta_K$ and $\norm{u_j}<r_K$, which implies $u_j=0$, a contradiction. Thus $\Stab_G(u_j)=H$ for all $j\ge j_0$, proving (i).

\emph{Step 4: Proof of the lower bound \eqref{eq:NH_ceiling}.}
When $|s_H|>B_H$, since $|s_H|$ and $B_H$ are integers, we have $|s_H|-B_H\ge1$. From \eqref{eq:count} we infer that $\sum_{t\in\mathcal{T}_H}\nu_H(t)\ge|s_H|-B_H\ge1$, and hence $\mathcal{T}_H\ne\varnothing$ and $M_H=\max_{t\in\Sigma_H\setminus\cE_H}\nu_H(t)\ge1$. Estimating each term in \eqref{eq:count} by $M_H$ yields
$$
 M_H\cdot\#\mathcal{T}_H\ge\sum_{t\in\mathcal{T}_H}\nu_H(t)\ge|s_H|-B_H\quad\Longrightarrow\quad
 N_H\ge\#\mathcal{T}_H\ge\frac{|s_H|-B_H}{M_H},
$$
from which \eqref{eq:NH_ceiling} follows since $N_H$ and $\#\mathcal{T}_H$ are integers.
\end{proof}

\begin{proof}[Proof of Theorem A]
Item (i) is a direct consequence of Theorem~\ref{thm:classical}(i) applied to $f_\lambda|_{\cH^H}$ together with Theorem~\ref{thm:PSC}, while Lemma~\ref{lem:uniqueness} rules out bifurcation at the endpoints $a$ and $b$. Item (ii) follows from Theorem~\ref{thm:budget}, and \eqref{eq:intro_morse} is simply \eqref{eq:sfmorse} for $L^H$.
\end{proof}

\begin{remark}\label{rem:edge_cases}
It is worth noticing that in \eqref{eq:BH_MH} we have $M_H=0$ if and only if $\Sigma_H\setminus\cE_H=\varnothing$; in that case $\Sigma_H=\cE_H$, and Lemma~\ref{lem:localflow} gives $|s_H|\le\sum_{t\in\Sigma_H}\nu_H(t)=B_H$, so that the condition $|s_H|>B_H$ automatically guarantees $M_H\ge1$. Notice also that since every subgroup $S>H$ contains some $K\in\mathcal K_H$, checking $\nu_K(t)=0$ for all $K\in\mathcal K_H$ is equivalent to checking $\nu_S(t)=0$ for all $S>H$. The example below shows that the estimate \eqref{eq:thmA_bound} is sharp even when $L_a$ and $L_b$ have non-trivial kernels on $\cH$.
\end{remark}

\begin{figure}[htbp]
\centering
\definecolor{gpwPath}{RGB}{24,52,73}
\definecolor{gpwBlue}{RGB}{33,103,181}
\definecolor{gpwTeal}{RGB}{22,133,128}
\definecolor{gpwOrange}{RGB}{195,106,22}
\resizebox{\linewidth}{!}{%
\begin{tikzpicture}[>=Latex,font=\small]
\path[use as bounding box] (0,-3.65) rectangle (14.8,6.95);
\node[text=gpwPath] at (7.4,6.7)
 {$G=C_4,\quad H=\langle r^2\rangle,\quad
   L_t^H=\operatorname{diag}(t+1,t-1,t)$};
\node at (3.2,6.05) {(a) A plane in the space of Hessians};
\node at (3.2,5.6)
 {$D(\alpha,\beta)=\operatorname{diag}(\alpha+1,\alpha-1,\beta)$};
\node at (10.4,6.05) {(b) Eigenvalues of the same path};

% The operator path in a two-dimensional affine plane of Hessians.
\begin{scope}[shift={(3.2,2.75)},x=1.15cm,y=1.15cm]
 \draw[step=1,gray!18,very thin] (-2,-2) grid (2,2);
 \draw[gray!50] (-2,-2) rectangle (2,2);
 \foreach \v in {-2,-1,0,1,2} {
  \draw[gray] (\v,-2)--(\v,-2.055);
  \node[below,text=gray] at (\v,-2.06) {$\v$};
  \draw[gray] (-2,\v)--(-2.055,\v);
  \node[left,text=gray] at (-2.06,\v) {$\v$};
 }
 \node at (0,-2.55) {$\alpha$};
 \node at (-2.63,0) {$\beta$};
 \draw[gpwBlue,dashed,thick] (-1,-2)--(-1,2);
 \draw[gpwBlue,dashed,thick] (1,-2)--(1,2);
 \draw[gpwOrange,dashed,thick] (-2,0)--(2,0);
 \draw[gpwPath,very thick] (-1.9,-1.9)--(1.9,1.9);
 \draw[gpwPath,very thick,->] (1.35,1.35)--(1.67,1.67);
 \filldraw[fill=gpwBlue,draw=white] (-1,-1) circle (2.6pt);
 \filldraw[fill=gpwOrange,draw=white] (0,0) circle (2.6pt);
 \filldraw[fill=gpwTeal,draw=white] (1,1) circle (2.6pt);
 \node[text=gpwBlue,fill=white,inner sep=1pt] at (-1.6,-.7) {$t=-1$};
 \node[text=gpwOrange,fill=white,inner sep=1pt] at (.45,-.32) {$t=0$};
 \node[text=gpwTeal,fill=white,inner sep=1pt] at (1.5,.72) {$t=1$};
\end{scope}

% Eigenvalues of L_t^H, with upward crossings at -1, 0 and 1.
\begin{scope}[shift={(10.4,2.75)},x=1.15cm,y=.72cm]
 \foreach \v in {-2,-1,0,1,2}
  {\draw[gray!18,very thin] (\v,-3.2)--(\v,3.2);}
 \foreach \v in {-3,-2,-1,0,1,2,3}
  {\draw[gray!18,very thin] (-2.05,\v)--(2.5,\v);}
 \draw[gray!50] (-2.05,-3.2) rectangle (2.5,3.2);
 \foreach \v in {-2,-1,0,1,2} {
  \draw[gray] (\v,-3.2)--(\v,-3.28);
  \node[below,text=gray] at (\v,-3.29) {$\v$};
 }
 \foreach \v in {-3,-2,-1,0,1,2,3} {
  \draw[gray] (-2.05,\v)--(-2.1,\v);
  \node[left,text=gray] at (-2.11,\v) {$\v$};
 }
 \node at (.2,-4.05) {parameter $t$};
 \node[rotate=90] at (-2.7,0) {eigenvalue};
 \draw[gpwPath] (-2.05,0)--(2.5,0);
 \draw[gpwBlue,very thick] (-2,-1)--(2,3);
 \draw[gpwOrange,very thick] (-2,-2)--(2,2);
 \draw[gpwTeal,very thick] (-2,-3)--(2,1);
 \draw[gpwBlue,dotted] (-1,-3.2)--(-1,0);
 \draw[gpwOrange,dotted] (0,-3.2)--(0,0);
 \draw[gpwTeal,dotted] (1,-3.2)--(1,0);
 \filldraw[fill=gpwBlue,draw=white] (-1,0) circle (2.6pt);
 \filldraw[fill=gpwOrange,draw=white] (0,0) circle (2.6pt);
 \filldraw[fill=gpwTeal,draw=white] (1,0) circle (2.6pt);
 \node[right,text=gpwBlue,inner sep=2pt] at (2,3) {$t+1$};
 \node[right,text=gpwOrange,inner sep=2pt] at (2,2) {$t$};
 \node[right,text=gpwTeal,inner sep=2pt] at (2,1) {$t-1$};
\end{scope}

\node[anchor=west] at (.35,-.65)
 {Dashed lines: $\det D(\alpha,\beta)=0$.};
\node[anchor=west,text=gpwOrange] at (.35,-1.05)
 {Orange line: $D|_{\mathcal H^G}=0$.};
\node[anchor=west,text=gpwPath] at (.35,-1.45)
 {Solid line: $L_t^H=D(t,t)$.};
\node[anchor=west,text=gpwBlue] at (7.6,-.65)
 {$t=-1:\ \ker L_{-1}^H=\mathbb R e_1$; exact $H$.};
\node[anchor=west,text=gpwOrange] at (7.6,-1.05)
 {$t=0:\ \ker L_0^H=\mathbb R e_y$; fixed by $G$.};
\node[anchor=west,text=gpwTeal] at (7.6,-1.45)
 {$t=1:\ \ker L_1^H=\mathbb R e_2$; exact $H$.};

% The singular and exceptional sets lie on the parameter interval.
\begin{scope}[shift={(7,-2.25)},x=2cm,y=1cm]
 \foreach \y in {0,-.75} {
  \draw[gray,->] (-2.1,\y)--(2.15,\y);
  \foreach \v in {-2,-1,0,1,2}
   {\draw[gray] (\v,{\y-.035})--(\v,{\y+.035});}
 }
 \node[left] at (-2.3,0) {$\Sigma_H$};
 \node[left] at (-2.3,-.75) {$\mathcal E_H$};
 \node[right] at (2.35,0) {$\{-1,0,1\}$};
 \node[right,text=gpwOrange] at (2.35,-.75) {$\{0\}$};
 \draw[gpwOrange,dotted] (0,0)--(0,-.75);
 \fill[gpwBlue] (-1,0) circle (2.6pt);
 \fill[gpwOrange] (0,0) circle (2.6pt);
 \fill[gpwTeal] (1,0) circle (2.6pt);
 \fill[gpwOrange] (0,-.75) circle (2.6pt);
 \foreach \v in {-2,-1,0,1,2}
  {\node[below,text=gray] at (\v,-.85) {$\v$};}
\end{scope}
\end{tikzpicture}%
}
\caption[Geometric interpretation of the singular and exceptional sets]{%
The restricted Hessian path in Example~\ref{ex:C4}.
In panel (a), the singular matrices in the affine plane
$D(\alpha,\beta)=\operatorname{diag}(\alpha+1,\alpha-1,\beta)$
form the lines $\alpha=\pm1$ and $\beta=0$.
The path $L_t^H=D(t,t)$ meets them at $t=-1,0,1$.
Panel (b) shows its eigenvalues; $e_1,e_2,e_y$ are the coordinate
vectors for $(x_1,x_2,y)$ on $\mathcal H^H$.
At $t=\pm1$, every nonzero kernel vector has exact stabilizer $H$;
at $t=0$, the kernel is fixed by $G$.
Thus $\Sigma_H=\{-1,0,1\}$ and $\mathcal E_H=\{0\}$.
Both are sets of parameter values. In general, membership in
$\mathcal E_H$ asserts the existence of a larger-symmetry kernel
direction, not that every kernel direction has larger symmetry.}
\label{fig:budget}
\end{figure}
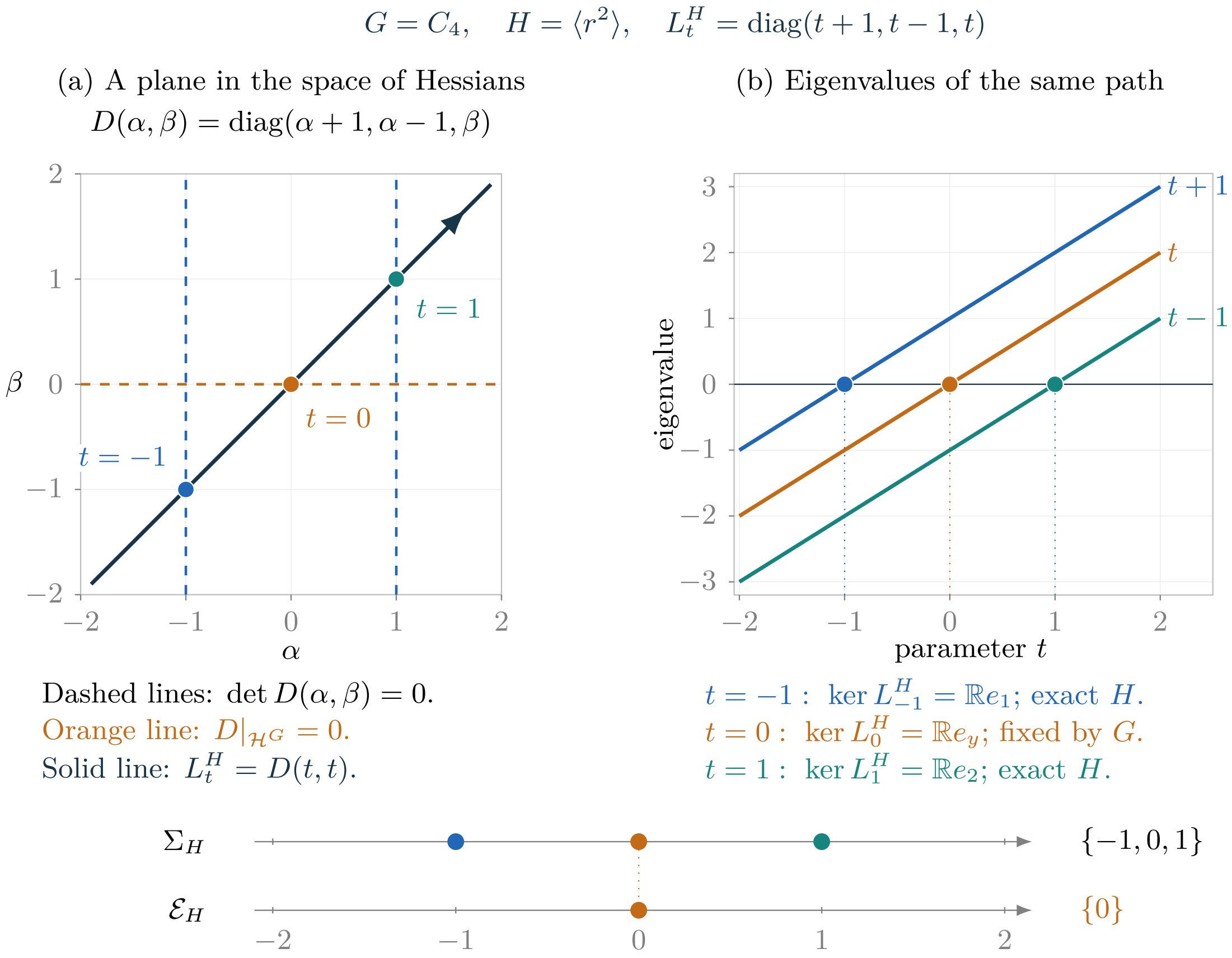

\begin{example}[Attaining the bound \eqref{eq:thmA_bound} with singular ambient endpoints]\label{ex:C4}
Let $G=C_4=\langle r\rangle$ act orthogonally on $\cH=\R^2\oplus\R\oplus\C\cong\R^5$ via
$$
 r(x_1,x_2,y,z)=(-x_1,-x_2,y,iz),
$$
and take $H=\langle r^2\rangle\cong C_2$. Then $\cH^H=\{(x_1,x_2,y,0)\}\cong\R^3$, $\cH^G=\{(0,0,y,0)\}\cong\R$, and $\mathcal{K}_H=\{G\}$. On $I=[-2,2]$, we consider the $G$-invariant family
$$
 f_t(x_1,x_2,y,z)=\frac{1}{2}\big((t+1)x_1^2+(t-1)x_2^2+ty^2+(4-t^2)|z|^2\big)
              +\frac{1}{4}\big(x_1^4+x_2^4+y^4+|z|^4\big).
$$
At $t=\pm2$, the kernel of $L_t$ on $\cH$ is the $2$-dimensional plane $\{0\}\times\{0\}\times\{0\}\times\C$, so $L_{\pm2}$ is degenerate. On $\cH^H$ and $\cH^G$, however, we have
$$
 L_t^H=\diag(t+1,t-1,t) \quad \textrm{ and } \quad L_t^G=t,
$$
which are both invertible at $t=\pm2$. As $t$ goes from $-2$ to $2$, all three eigenvalues of $L_t^H$ cross zero from negative to positive, so that
$$
 s_H=m^-(L_{-2}^H)-m^-(L_2^H)=3,\qquad
 \Sigma_H=\{-1,0,1\},\qquad \cE_H=\{0\},\qquad B_H=\nu_H(0)=1,\qquad M_H=1,
$$
and Theorem A(ii) gives $N_H\ge\lceil(3-1)/1\rceil=2$. Solving $\nabla f_t=0$ explicitly shows that near $t=-1$, the non-trivial critical points are $(\pm\sqrt{-t-1},0,0,0)$ with exact stabilizer $H$; near $t=1$, they are $(0,\pm\sqrt{1-t},0,0)$ with exact stabilizer $H$; whereas near the exceptional parameter $t=0$, for sufficiently small solutions the equations $(t+1)x_1+x_1^3=0$, $(t-1)x_2+x_2^3=0$, and $(4-t^2+|z|^2)z=0$ force $x_1=x_2=z=0$, so the only bifurcating critical points are $(0,0,\pm\sqrt{-t},0)\in\cH^G$ with stabilizer $G\ne H$. Thus the exceptional crossing at $t=0$ produces no exact-$H$ solutions, and the lower bound $N_H=2$ is attained.
\end{example}
Figure~\ref{fig:budget} gives a geometric interpretation of the
sets $\Sigma_H$ and $\cE_H$ in Example~\ref{ex:C4}, showing both
the path of restricted Hessians and its eigenvalues.

% Paste this entire block in the main manuscript, replacing the figure
% labelled fig:budget. The existing preamble already loads all packages.
% No external image, PDF, or additional source file is needed.

\subsection{Equivariant reduction and simple branches}\label{subsec:reduction}

Theorem A yields sequences of critical points accumulating at $(t,0)$ under $C^0$ dependence in the parameter $\lambda$. If one wants to get smooth local curves of critical points, one needs higher differentiability in $\lambda$ together with a $1$-dimensional transverse kernel. With a view to the applications in Sections 6 and 7, we briefly recall here the standard equivariant Lyapunov--Schmidt reduction (cf.\ \cite[Lemma 5.5]{IJSW}) and its Crandall--Rabinowitz consequence \cite[Theorem 1.7]{CrandallRabinowitz}.

\begin{theorem}[Equivariant Lyapunov--Schmidt reduction]\label{thm:LS}
Under \textup{(H1)--(H3)}, fix $t_*\in(a,b)$, set $E=\ker L_{t_*}$, and let $P\in\mathcal{L}(\cH)$ be the orthogonal projection onto $E$ and $Q=I-P$. Then there exist open neighborhoods $J_0\subset I$ of $t_*$, $U_0\subset E$ of $0$, and $W_0\subset E^\perp$ of $0$, and a unique map $w:J_0\times U_0\to W_0$, continuous in $\lambda$ and $C^1$ in $z\in U_0$ with continuous derivative $D_zw$, such that
$$
 Q\nabla f_\lambda(z+w(\lambda,z))=0,\qquad w(\lambda,0)=0,\qquad D_zw(t_*,0)=0,
$$
and all critical points of $f_\lambda$ in $U_0\oplus W_0$ are of the form $u=z+w(\lambda,z)$. Moreover, the reduced functional $r_\lambda(z):=f_\lambda(z+w(\lambda,z))-f_\lambda(0)$ on $U_0$ is $G$-invariant, $C^2$ in $z$, and satisfies
$$
 \nabla r_\lambda(z)=P\nabla f_\lambda(z+w(\lambda,z))
 \quad \textrm{ and } \quad
 \Stab_G(z+w(\lambda,z))=\Stab_G(z).
$$
If in addition $(\lambda,u)\mapsto\nabla f_\lambda(u)$ is of class $C^k$ ($k\ge1$), then $w$, $r_\lambda$, and $\nabla r_\lambda$ are of class $C^k$, and $z\in U_0\cap E^H$ implies $w(\lambda,z)\in(E^\perp)^H$ for any $H\le G$.
\end{theorem}
\begin{proof}
Since $L_{t_*}$ is self-adjoint and Fredholm, $QL_{t_*}|_{E^\perp}:E^\perp\to E^\perp$ is a toplinear isomorphism, and hence the implicit function theorem applied to $(\lambda,z,w)\mapsto Q\nabla f_\lambda(z+w)$ at $(t_*,0,0)$ gives the unique local solution $w(\lambda,z)$ with $D_zw(t_*,0)=0$. Because $P$ and $Q$ commute with the action of $G$, uniqueness implies that $w(\lambda,gz)=gw(\lambda,z)$, $\forall\,g\in G$. Since $E$ and $E^\perp$ are orthogonal and $G$-invariant, $g(z+w(\lambda,z))=z+w(\lambda,z)$ holds if and only if $gz=z$ and $gw(\lambda,z)=w(\lambda,z)$, which proves that $\Stab_G(z+w(\lambda,z))=\Stab_G(z)$ and $w(\lambda,U_0\cap E^H)\subset(E^\perp)^H$. Finally, differentiating $r_\lambda$ with respect to $z$ gives
$$
 Dr_\lambda(z)[v]=\ip{\nabla f_\lambda(z+w(\lambda,z))}{v+D_zw(\lambda,z)v}
 =\ip{P\nabla f_\lambda(z+w(\lambda,z))}{v},
$$
since $Q\nabla f_\lambda(z+w(\lambda,z))=0$ and $\ran D_zw(\lambda,z)\subset E^\perp$. Differentiating $\nabla r_\lambda(z)=P\nabla f_\lambda(z+w(\lambda,z))$ once more in $z$ shows that $r_\lambda$ is $C^2$ in $z$, and when $(\lambda,u)\mapsto\nabla f_\lambda(u)$ is $C^k$, the identity $r_\lambda(z)=\int_0^1\ip{\nabla r_\lambda(tz)}{z}\,dt$ (which follows from $r_\lambda(0)=0$) implies that $r_\lambda$ is $C^k$ as well.
\end{proof}

Whenever $\ker L_{t_*}^H$ is $1$-dimensional and the crossing form is non-zero, Crandall--Rabinowitz's theorem \cite[Theorem 1.7]{CrandallRabinowitz} yields a smooth local curve of critical points with exact stabilizer $H$ (cf.\ also \cite{Cicogna1981,Vanderbauwhede1982,GSS1988}).

\begin{corollary}[Simple fixed-space crossing]\label{cor:simple}
Under \textup{(H1)--(H3)}, suppose that $F(\lambda,u)=\nabla f_\lambda(u)$ be of class $C^k$ ($k\ge2$) near $(t_*,0)$ for some $t_*\in(a,b)$, and that for a subgroup $H\le G$ we have
$$
 \ker L_{t_*}^H=\R v,\qquad \norm{v}=1,\qquad
 \ip{\partial_\lambda L_\lambda|_{t_*}v}{v}\ne0,\qquad \Stab_G(v)=H.
$$
Then there exist $\varepsilon_0>0$ and a $C^{k-1}$ curve $(-\varepsilon_0,\varepsilon_0)\ni s\mapsto(\lambda(s),u(s))\in I\times\cH^H$ such that $\lambda(0)=t_*$, $u(s)=sv+o(s)$ as $s\to0$, $\nabla f_{\lambda(s)}(u(s))=0$, and $\Stab_G(u(s))=H$ for all $0<|s|<\varepsilon_0$. Moreover, near $(t_*,0)$, all non-trivial critical points of $f_\lambda$ in $\cH^H$ belong to this curve.
\end{corollary}
\begin{proof}
On the subspace $\cH^H$ we have $\ran L_{t_*}^H=(\R v)^\perp\cap\cH^H$, so that $\ip{\partial_\lambda L_\lambda|_{t_*}v}{v}\ne0$ is equivalent to the transversality condition $\partial_\lambda D_uF(t_*,0)v\notin\ran L_{t_*}^H$ of \cite[Theorem 1.7]{CrandallRabinowitz}. By Theorem~\ref{thm:PSC}, the resulting $C^{k-1}$ branch $(\lambda(s),u(s))$ in $\cH^H$ consists of critical points of $f_\lambda$ on $\cH$. Finally, if $H\ne G$, for every $g\in G\setminus H$ we have $d_g:=\norm{gv-v}>0$; setting $d_*:=\min_{g\in G\setminus H}d_g>0$, since $u(s)/s\to v$ as $s\to0$, for $0<|s|<\varepsilon_0$ small enough we get
$$
 \norm{g\frac{u(s)}{s}-\frac{u(s)}{s}}
 \ge d_g-2\norm{\frac{u(s)}{s}-v}\ge\frac{d_*}{2}>0, \quad \forall\, g\in G\setminus H,
$$
and hence $\Stab_G(u(s))=H$.
\end{proof}

\subsection{Comparison with representation-valued spectral flow and equivariant degree}\label{subsec:limitations}

We close this section with a few remarks on the assumptions of Theorem A and on its relation to $RO(G)$-valued spectral flow.

First of all, let us point out that in Theorem A the ambient operators $L_a,L_b$ on $\cH$ may be degenerate, provided that their restrictions $L_a^H,L_b^H$ to $\cH^H$ be invertible. For instance, if $G=C_2$ acts on $\R^2$ by $(x,y)\mapsto(x,-y)$ and $f_\lambda(x,y)=\frac{1}{2}\lambda x^2+\frac{1}{4}x^4+\frac{1}{4}y^4$ on $[-1,1]$, then $L_\lambda=\diag(\lambda,0)$ has non-trivial kernel for all $\lambda\in[-1,1]$, whereas $L_\lambda^G=\lambda$ on $\cH^G=\{y=0\}$ is invertible at $\lambda=\pm1$ with $s_G=1$ and detects the branch $(\pm\sqrt{-\lambda},0)$. On the contrary, the assumption that $L_a^H,L_b^H$ be invertible on $\cH^H$ cannot be dropped, as shown by $f_\lambda(x)=\frac{1}{4}x^4-\frac{1}{2}\lambda x^2$ on $[0,1]$, where $m^-(L_0)=0$ and $m^-(L_1)=1$, but bifurcation occurs at the endpoint $\lambda=0$ and not in $(0,1)$.

Second, it is natural to ask whether a non-zero spectral flow on the orthogonal complement $(\cH^G)^\perp$ of the fixed-point space $\cH^G$ suffices by itself to produce symmetry breaking when $L_{t_*}^G$ is degenerate. As the next example shows, this is not true in general.

\begin{example}[Transverse index jump without symmetry breaking]\label{ex:transverse_no_sb}
Consider $G=C_2$ acting on $\R\oplus\R^r$ via $(x,v)\mapsto(x,-v)$, let $c>0$, and define on $[-1,1]$ the $G$-invariant polynomial functional
$$
 f_\lambda(x,v)=\frac{1}{3}x^3-\frac{1}{2}\lambda x^2+x|v|^2+\frac{1}{2}c\lambda|v|^2.
$$
At $u=(0,0)$ we have $L_\lambda=\diag(-\lambda,c\lambda I_r)$, whose restriction $c\lambda I_r$ to $\{0\}\oplus\R^r$ has spectral flow $r\ne0$ across $\lambda=0$. Nevertheless, the critical-point equations read $x^2-\lambda x+|v|^2=0$ and $(2x+c\lambda)v=0$. If $v\ne0$, then $x=-c\lambda/2$, and substituting into the first equation gives $\frac{1}{4}c(c+2)\lambda^2+|v|^2=0$, which has no solution with $v\ne0$ as $c>0$. Hence all critical points have $v=0$ (namely $(0,0)$ and $(\lambda,0)$) and stabilizer $G$. Notice that for $H=\{e\}$ in Theorem A we have $\Sigma_{\{e\}}=\{0\}=\cE_{\{e\}}$, $s_{\{e\}}=r-1$, and $B_{\{e\}}=\nu_{\{e\}}(0)=r+1\ge|s_{\{e\}}|$, so that the condition $|s_{\{e\}}|>B_{\{e\}}$ is not satisfied.
\end{example}

\begin{remark}\label{rem:comparison}
When $L_a$ and $L_b$ are invertible on $\cH$, the $G$-equivariant spectral flow $\sflo_G(L;[a,b])=[U]-[V]\in RO(G)$ introduced in \cite{IJW2021,IJSW,IJSW2026} is related to $s_H$ through the homomorphism $\Phi_H([U]-[V])=\dim U^H-\dim V^H$, namely
$$
 \Phi_H\big(\sflo_G(L;[a,b])\big)=\sflo(L^H;[a,b])=s_H.
$$
As already observed in the introduction, for cyclic groups such as $G=C_5$ the two non-isomorphic $2$-dimensional rotation representations $E_1$ (angle $2\pi/5$) and $E_2$ (angle $4\pi/5$) have $\dim E_1^K=\dim E_2^K$ for both $K=\{e\}$ and $K=C_5$, so that $L_t=tI_{E_1}\oplus(-tI_{E_2})$ on $[-1,1]$ has $\sflo_{C_5}(L;[-1,1])=[E_1]-[E_2]\ne0$ in $RO(C_5)$ while $s_K=0$ for all $K\le C_5$. Conversely, for groups such as $G=D_4$, decomposing any real $D_4$-module $V$ into the four one-dimensional characters $A_1,A_2,B_1,B_2$ (where $(r,s)$ acts by $(+,+),(+,-),(-,+),(-,-)$) and the two-dimensional representation $E$ shows that the fixed-space dimensions $d_K=\dim V^K$ uniquely determine all multiplicities via
\[
\begin{split}
 m_{A_1}&=d_{D_4},\qquad m_{A_2}=d_{\langle r\rangle}-d_{D_4},\\
 m_{B_1}&=d_{\langle r^2,s\rangle}-d_{D_4},\qquad
 m_{B_2}=d_{\langle r^2,rs\rangle}-d_{D_4},\\
 2m_E&=\dim V-m_{A_1}-m_{A_2}-m_{B_1}-m_{B_2},
\end{split}
\]
which implies in particular that $D_4$ satisfies the ``nice group'' hypothesis of \cite[Theorem 3.1]{IJSW}. Within this setting, Theorem A complements both $RO(G)$-valued spectral flow \cite{IJSW} and global equivariant degree methods \cite{IzeVignoliBook,GarciaAzpeitiaIze2011,GarciaAzpeitiaIze,Kowalczyk} by providing an explicit lower bound \eqref{eq:thmA_bound} on the number of exact-$H$ bifurcation instants under restricted endpoint invertibility $\ker L_a^H=\ker L_b^H=\{0\}$.
\end{remark}

\section{Variational coordinates for varying form domains}\label{sec:form}

Theorem A is stated for $C^2$ functionals on a fixed Hilbert space $\cH$. In several variational problems associated with unbounded self-adjoint operators $A_\lambda$ on $\cH$, the natural domain of the energy functional is the form domain $V_\lambda=D(|A_\lambda|^{1/2})$, which may depend on $\lambda$. In this section we show how the canonical isometries $U_\lambda:\cH\to V_\lambda$ allow us to bring the problem back to the fixed space $\cH$ (cf.\ \cite{CandelaWaterstraat,Pejsachowicz2025}).

\subsection{Spectral forms and operator topologies}\label{subsec:forms}

Let $A:D(A)\subset\cH\to\cH$ be a densely defined self-adjoint operator (not necessarily semibounded) on a real separable Hilbert space $\cH$. Using the Borel functional calculus, we define the form domain $V_A=D(|A|^{1/2})$ and equip it with the two equivalent Hilbert norms
\[
 \norm{u}_{V_A}^2=\norm{u}^2+\norm{|A|^{1/2}u}^2,\qquad
 \norm{u}_{A,*}=\norm{(1+A^2)^{1/4}u}.
\]
The associated symmetric bilinear spectral form on $V_A$ is given by
\[
 q_A(u,v)=\ip{\sgn(A)|A|^{1/2}u}{|A|^{1/2}v}\qquad(u,v\in V_A),
\]
where $\sgn(t)=t/|t|$ for $t\ne0$ and $\sgn(0)=0$.

\begin{lemma}[Canonical form coordinates]\label{lem:form}
For any densely defined self-adjoint operator $A$ on $\cH$, the following properties hold:
\begin{enumerate}[(i)]
\item $V_A$ is a Hilbert space under both $\norm{\cdot}_{V_A}$ and $\norm{\cdot}_{A,*}$, and
\[
 \norm{u}_{A,*}\le\norm{u}_{V_A}\le 2^{1/4}\norm{u}_{A,*}\qquad\text{for all }u\in V_A.
\]
\item The canonical transport operator $U_A=(1+A^2)^{-1/4}:\cH\to V_A$ is a surjective isometry with respect to $\norm{\cdot}_{A,*}$ (and a topological isomorphism with respect to $\norm{\cdot}_{V_A}$) that pulls $q_A$ back to the bounded Riesz transform $\Psi(A)=A(1+A^2)^{-1/2}\in\mathcal{L}(\cH)$:
\begin{equation}\label{eq:transport}
 q_A(U_Ax,U_Ay)=\ip{\Psi(A)x}{y}\qquad\text{for all }x,y\in\cH.
\end{equation}
\item The operators $A$ and $\Psi(A)$ have the exact same kernel and negative spectral subspace, and $A$ is Fredholm if and only if $\Psi(A)$ is Fredholm.
\item For any $u\in V_A$ and $f\in\cH$,
\begin{equation}\label{eq:weakstrong}
 q_A(u,v)=\ip{f}{v}, \quad \forall\, v\in V_A
 \quad\Longleftrightarrow\quad u\in D(A)\text{ and }Au=f.
\end{equation}
\end{enumerate}
\end{lemma}
\begin{proof}
Denoting by $E_A$ the spectral measure of $A$ and setting $\mu_u(B)=\norm{E_A(B)u}^2$, we have $u\in V_A$ if and only if $\int_{\R}|t|\,d\mu_u(t)<\infty$, and $\norm{u}_{A,*}^2=\int_{\R}\sqrt{1+t^2}\,d\mu_u(t)$. Item (i) follows at once from the elementary inequalities $\sqrt{1+t^2}\le1+|t|\le\sqrt{2}\sqrt{1+t^2}$ and the fact that $|A|^{1/2}$ is closed. Moreover, for any $x\in\cH$ we compute $\norm{U_Ax}_{A,*}^2=\int_{\R}\sqrt{1+t^2}(1+t^2)^{-1/2}\,d\mu_x(t)=\norm{x}^2$, while any $u\in V_A$ can be written as $u=U_Ax$ with $x=(1+A^2)^{1/4}u\in\cH$. Since $\sgn(t)|t|^{1/2}(1+t^2)^{-1/4}\cdot|t|^{1/2}(1+t^2)^{-1/4}=t(1+t^2)^{-1/2}=\Psi(t)$, we obtain \eqref{eq:transport}. As $t\mapsto\Psi(t)=t/\sqrt{1+t^2}$ is a strictly increasing homeomorphism of $\R$ onto $(-1,1)$ with $\Psi(0)=0$, $A$ and $\Psi(A)$ have the same kernel, the same negative spectral subspace $\ran\mathbf{1}_{(-\infty,0)}(A)=\ran\mathbf{1}_{(-1,0)}(\Psi(A))$, and the same Fredholm property. Finally, if $v\in D(A)$ then $q_A(u,v)=\ip{u}{Av}$; thus $q_A(u,v)=\ip{f}{v}$ for all $v\in D(A)$ means that $u\in D(A^*)=D(A)$ and $Au=f$, and the density of $D(A)$ in $V_A$ gives \eqref{eq:weakstrong} (cf.\ \cite{Kato}).
\end{proof}

For unbounded self-adjoint operators continuity in the gap metric is equivalent to norm-resolvent continuity $\lambda\mapsto(A_\lambda-iI)^{-1}$, whereas Riesz continuity (namely, norm continuity of $\lambda\mapsto\Psi(A_\lambda)$) is strictly stronger in general \cite{BLP,Lesch}. Nevertheless, the two topologies coincide on uniformly semibounded families.

\begin{lemma}[Uniform semiboundedness implies Riesz continuity]\label{lem:semibounded}
Let $\lambda\mapsto A_\lambda$ be a gap-continuous path of densely defined self-adjoint operators on $\cH$. Then $\lambda\mapsto f(A_\lambda)$ is norm-continuous in $\mathcal{L}(\cH)$ for every $f\in C_0(\R)$ and every continuous function $f$ having equal finite limits at $\pm\infty$. In particular, $\lambda\mapsto U_{A_\lambda}=(1+A_\lambda^2)^{-1/4}$ is norm-continuous in $\mathcal{L}(\cH)$. Furthermore, if there is $C>0$ such that $A_\lambda\ge-CI$ for all $\lambda\in I$ (resp.\ $A_\lambda\le CI$ for all $\lambda\in I$), then $\lambda\mapsto\Psi(A_\lambda)$ is norm-continuous on $I$.
\end{lemma}
\begin{proof}
Since $\lambda\mapsto(A_\lambda\pm iI)^{-1}$ is norm-continuous, the same holds for $p((A_\lambda-iI)^{-1},(A_\lambda+iI)^{-1})$ for any polynomial $p$. By Stone--Weierstrass's theorem, such polynomials are uniformly dense in $C_0(\R)$, and hence the inequality $\norm{f(A_\lambda)-h(A_\lambda)}\le\norm{f-h}_\infty$ implies that $\lambda\mapsto f(A_\lambda)$ is norm-continuous for all $f\in C_0(\R)\oplus\R 1$, which includes $(1+t^2)^{-1/4}$. If now $\spec(A_\lambda)\subset[-C,\infty)$ for all $\lambda\in I$, it suffices to pick a continuous function $h\in C_0(\R)\oplus\R 1$ such that $h(t)=t/\sqrt{1+t^2}$ for $t\ge-C$ and $\lim_{t\to\pm\infty}h(t)=1$: then $\Psi(A_\lambda)=h(A_\lambda)$ is norm-continuous on $I$. The case in which $A_\lambda\le CI$ is proved in the same way.
\end{proof}

Whenever $A_\lambda$ is $G$-equivariant, $P_H$ preserves $D(A_\lambda)$ and commutes with the spectral measure of $A_\lambda$, so that $A_\lambda^H=A_\lambda|_{D(A_\lambda)\cap\cH^H}$ is self-adjoint on $\cH^H$ and satisfies
$$
 \Psi(A_\lambda^H)=\Psi(A_\lambda)\big|_{\cH^H},\qquad
 U_{A_\lambda^H}=U_{A_\lambda}\big|_{\cH^H},\qquad
 \ker A_\lambda^H=(\ker A_\lambda)^H,
$$
and for a Riesz-continuous Fredholm path we set
\begin{equation}\label{eq:sftransform}
 \sflo(A^H;[a,b]):=\sflo\big(\Psi(A^H);[a,b]\big).
\end{equation}

\subsection{Transport theorem for varying form domains}\label{subsec:transport}

We consider a $G$-invariant family of functionals
$$
 \Phi_\lambda(u)=\frac{1}{2}q_{A_\lambda}(u,u)+N_\lambda(u),\qquad u\in U_\lambda(B_\rho(0))\subset V_\lambda=D(|A_\lambda|^{1/2}),
$$
under the following hypotheses:
\begin{enumerate}[(V1)]
\item For every $\lambda\in I=[a,b]$, $A_\lambda:D(A_\lambda)\subset\cH\to\cH$ is a densely defined, $G$-equivariant, self-adjoint Fredholm operator, and $\lambda\mapsto\Psi(A_\lambda)$ is norm-continuous in $\mathcal{L}(\cH)$.
\item For each $\lambda\in I$, $N_\lambda$ is a $G$-invariant functional defined on $U_\lambda(B_\rho(0))\subset V_\lambda$, where $U_\lambda=(1+A_\lambda^2)^{-1/4}$ and $B_\rho(0)\subset\cH$, such that the pulled-back family $\widehat{N}_\lambda:=N_\lambda\circ U_\lambda:B_\rho(0)\to\R$ satisfies hypothesis (H1).
\item For all $\lambda\in I$, $D\widehat{N}_\lambda(0)=0$ and $D^2\widehat{N}_\lambda(0)=0$.
\end{enumerate}

\begin{theorem}[Fixed-space bifurcation via canonical form transport]\label{thm:transport}
Assume hypotheses \textup{(V1)--(V3)} hold, fix $H\le G$, and suppose that $A_a^H$ and $A_b^H$ have trivial kernel.
\begin{enumerate}[(i)]
\item If $s_H=\sflo(A^H;[a,b])\ne0$, then there exist $t\in(a,b)$, $\lambda_j\to t$, and non-zero $H$-fixed form-critical points $u_j\in V_{\lambda_j}\cap\cH^H\setminus\{0\}$ of $\Phi_{\lambda_j}$ (meaning $\tfrac{1}{2}Dq_{A_{\lambda_j}}(u_j)[v]+DN_{\lambda_j}(u_j)[v]=0$ for all $v\in V_{\lambda_j}$) such that $\norm{U_{\lambda_j}^{-1}u_j}_{\cH}\to0$ and $\norm{u_j}_{\cH}\to0$.
\item If $\Sigma_H=\{t\in(a,b):\dim\ker A_t^H>0\}$ is finite and $|s_H|>B_H$, where $\cE_H,B_H,M_H$ are defined by \eqref{eq:exceptional}--\eqref{eq:BH_MH} using $\nu_K(t)=\dim\ker A_t^K$, then at least $\lceil(|s_H|-B_H)/M_H\rceil$ distinct parameters $t\in\Sigma_H\setminus\cE_H$ admit such sequences with exact stabilizer $\Stab_G(u_j)=H$.
\item If in addition $DN_\lambda(u)[v]=\ip{F_\lambda(u)}{v}_{\cH}$ for some map $F_\lambda:U_\lambda(B_\rho(0))\to\cH$, then every form-critical point $u_j$ belongs to the operator domain $D(A_{\lambda_j})$ and solves the strong equation $A_{\lambda_j}u_j+F_{\lambda_j}(u_j)=0$.
\end{enumerate}
\end{theorem}
\begin{proof}
Pulling $\Phi_\lambda$ back to $B_\rho(0)\subset\cH$ via $U_\lambda$ defines the functional
\[
 f_\lambda(x)=\Phi_\lambda(U_\lambda x)=\frac{1}{2}\ip{\Psi(A_\lambda)x}{x}+\widehat{N}_\lambda(x).
\]
In view of Lemma~\ref{lem:form} and (V1)--(V3), the family $f_\lambda$ fulfills (H1)--(H3) on $\cH$ with $L_\lambda=\Psi(A_\lambda)$ and $\dim\ker L_\lambda^K=\dim\ker A_\lambda^K=\nu_K(\lambda)$. Therefore, Theorem A applied to $f_\lambda$ gives critical points $x_j\in B_\rho(0)\cap\cH^H\setminus\{0\}$ of $f_{\lambda_j}$ with $x_j\to0$ and $\lambda_j\to t\in(a,b)$ having the required stabilizer and multiplicity properties. If we set $u_j=U_{\lambda_j}x_j\in V_{\lambda_j}\cap\cH^H\setminus\{0\}$, then the fact that $U_{\lambda_j}:\cH\to V_{\lambda_j}$ is a $G$-equivariant isomorphism implies $D\Phi_{\lambda_j}(u_j)[U_{\lambda_j}y]=Df_{\lambda_j}(x_j)[y]=0$ for all $y\in\cH$, $\Stab_G(u_j)=\Stab_G(x_j)$, and $\norm{u_j}_{\cH}\le\norm{U_{\lambda_j}^{-1}u_j}_{\cH}=\norm{x_j}_{\cH}\to0$. Item (iii) is then an immediate consequence of \eqref{eq:weakstrong}.
\end{proof}

\section{Regular Sturm--Liouville systems and boundary paths}\label{sec:SL}

In the case of regular differential systems, one can usually work directly on a Sobolev space without having to use the fractional powers of Section 3. Here we deduce explicit criteria for matrix Sturm--Liouville problems, recall how moving boundary conditions of constant rank can be trivialized in $H^1$, and examine a Robin-to-Dirichlet path.

\subsection{Coefficient criteria and matrix Morse indices}\label{subsec:concrete}

We set $V_D=H^1_0(0,1;\R^d)$ and $V_R=\{u\in H^1(0,1;\R^d):u(1)=0\}$, endowed with the norm $\norm{u}_V=\norm{u'}_{L^2}$, and consider a potential $W:I\times[0,1]\times\R^d\to\R$ such that $W$, $\nabla_uW$, and $D_u^2W$ be continuous on $I\times[0,1]\times\R^d$ with
\begin{equation}\label{eq:matrixnonlinearity}
 \nabla_uW(\lambda,t,0)=0 \quad \textrm{ and } \quad D_u^2W(\lambda,t,0)=0, \quad \forall\, (\lambda,t)\in I\times[0,1].
\end{equation}

\begin{corollary}[Regular matrix Sturm--Liouville systems]\label{cor:matrixregular}
Let $P_\lambda(t),R_\lambda(t)$ be symmetric $d\times d$ matrices on $[0,1]$ such that $\lambda\mapsto P_\lambda$ be continuous in $C^1([0,1];\R^{d\times d})$ with $P_\lambda(t)\ge p_*I_d>0$, and $(\lambda,t)\mapsto R_\lambda(t)$ be continuous. Impose $u(1)=0$ at $t=1$ and either the Dirichlet condition $u(0)=0$ or the Robin condition $P_\lambda(0)u'(0)=B_\lambda u(0)$ at $t=0$, where $\lambda\mapsto B_\lambda$ is a continuous path of symmetric $d\times d$ matrices. Denoting by $A_\lambda$ the self-adjoint realization of $-(P_\lambda u')'+R_\lambda u$ in $L^2(0,1;\R^d)$ with these boundary conditions, and setting $m_\lambda=m^-(A_\lambda)$ and $\nu_\lambda=\dim\ker A_\lambda$, the following hold:
\begin{enumerate}[(i)]
\item If $A_a,A_b$ are invertible and $m_a\ne m_b$, then non-trivial solutions of
$$
 -(P_\lambda u')'+R_\lambda u+\nabla_uW(\lambda,t,u)=0
$$
bifurcate from $u=0$ in $H^1(0,1;\R^d)$ at some instant $t\in(a,b)$.
\item Under the assumptions of item \textup{(i)}, if $\Sigma=\{\lambda\in(a,b):\nu_\lambda>0\}$ is finite, then $\Sigma\ne\varnothing$, $M=\max_{\lambda\in\Sigma}\nu_\lambda\le d$, and the number of distinct bifurcation instants satisfies
\begin{equation}\label{eq:matrixcount}
 N_{\mathrm{bif}}\ge\ceil{\frac{|m_a-m_b|}{M}}\ge\ceil{\frac{|m_a-m_b|}{d}}.
\end{equation}
\item In the Dirichlet case $u(0)=u(1)=0$, if $Y_\lambda(t)\in\R^{d\times d}$ is the matrix solution of $-(P_\lambda Y_\lambda')'+R_\lambda Y_\lambda=0$ with $Y_\lambda(0)=0$ and $P_\lambda(0)Y_\lambda'(0)=I_d$, then
\begin{equation}\label{eq:matrixshooting}
 m_\lambda=\sum_{t\in(0,1)}\dim\ker Y_\lambda(t) \quad \textrm{ and } \quad \nu_\lambda=\dim\ker Y_\lambda(1).
\end{equation}
\end{enumerate}
\end{corollary}
\begin{proof}
On $V=V_D$ resp.\ $V_R$, we define the action functional
$$
 f_\lambda(u)=\frac{1}{2}q_\lambda(u,u)+\int_0^1W(\lambda,t,u(t))\,dt,
$$
where $q_\lambda(u,v)=\int_0^1(\ip{P_\lambda u'}{v'}+\ip{R_\lambda u}{v})\,dt+\ip{B_\lambda u(0)}{v(0)}$ in the Robin case (the boundary term being omitted when $V=V_D$). Since the embedding $V\hookrightarrow C^0([0,1];\R^d)$ is compact, a standard computation shows that $f_\lambda$ fulfills (H1)--(H3) on $V$ and that its Riesz Hessian $L_\lambda\in\mathcal{L}(V)$ is a compact perturbation of a positive isomorphism, with $m^-(L_\lambda)=m^-(A_\lambda)=m_\lambda$ and $\dim\ker L_\lambda=\dim\ker A_\lambda=\nu_\lambda\le d$ \cite{Kato,Zettl}. Items (i)--(ii) then follow from Theorem~\ref{thm:classical}, whereas \eqref{eq:matrixshooting} is the well-known Morse oscillation theorem for matrix Sturm--Liouville systems (cf.\ \cite{BOPW21,HowardOscillation,Zettl}).
\end{proof}

Notice that when $d=1$ we have $M=1$, and by Sturm's theorem $m_\lambda$ is simply the number of interior nodes in $(0,1)$ of the solution $y_\lambda$ satisfying $(y_\lambda(0),p_\lambda(0)y_\lambda'(0))=(0,1)$ resp.\ $(1,B_\lambda)$. Likewise, if $P$ and the boundary conditions do not depend on $\lambda$ and $R_\lambda(t)=R_0(t)-\lambda D(t)$ with $D(t)\ge d_*I_d>0$, then the crossing form $\Gamma(v,v)=-\int_0^1\ip{Dv}{v}\,dt$ is negative definite on the kernel, so that every eigenvalue $\lambda_n\in(a,b)$ of $A_0u=\lambda_n Du$ has non-zero local spectral flow and is a bifurcation instant.

\subsection{Moving essential boundary constraints}\label{subsec:regularforms}

If the essential boundary conditions depend on $\lambda$ but have constant rank, one can trivialize the family of energy spaces in $H^1$ by means of Kato's transformation \cite{Kato}. More precisely, let $X=H^1(0,1;\R^d)$, let $\gamma u=(u(0),u(1))\in\R^{2d}$ be the trace map, and consider boundary conditions of the form
$$
 B_\lambda u:=C_\lambda\gamma u=0,
$$
where $k$ is fixed and $\lambda\mapsto C_\lambda\in\mathcal L(\R^{2d},\R^k)$ is a norm-continuous path of surjective matrices. Endowing $X$ with a $G$-invariant $H^1$ scalar product, the orthogonal projection onto $V_\lambda=\ker B_\lambda\subset X$ is given by
$$
 P_\lambda=I-B_\lambda^*(B_\lambda B_\lambda^*)^{-1}B_\lambda\in\mathcal{L}(X),
$$
and depends norm-continuously on $\lambda$. On any subinterval on which $\norm{P_\lambda-P_{\lambda_0}}<1$, Kato's operator \cite{Kato}
$$
 T_\lambda=P_\lambda P_{\lambda_0}+(I-P_\lambda)(I-P_{\lambda_0})=I+(P_\lambda-P_{\lambda_0})(2P_{\lambda_0}-I)\in\mathcal{L}(X)
$$
is invertible and satisfies $P_\lambda T_\lambda=T_\lambda P_{\lambda_0}$. Thus $S_\lambda=T_\lambda|_{V_{\lambda_0}}:V_{\lambda_0}\to V_\lambda$ is a linear isomorphism depending norm-continuously on $\lambda$ and commuting with $G$ whenever $B_\lambda$ is $G$-equivariant; patching such maps over a finite cover of $I$ identifies all subspaces $V_\lambda$ with $V_a\subset H^1\hookrightarrow C^0$.

\subsection{Robin paths and computable symmetry breaking}\label{subsec:robin}

The constant-rank assumption of Subsection~\ref{subsec:regularforms} is violated when a Robin boundary condition degenerates into the Dirichlet condition. The next proposition shows that along such a path gap continuity does not imply Riesz continuity, and the Morse index can jump even though the operators remain invertible.

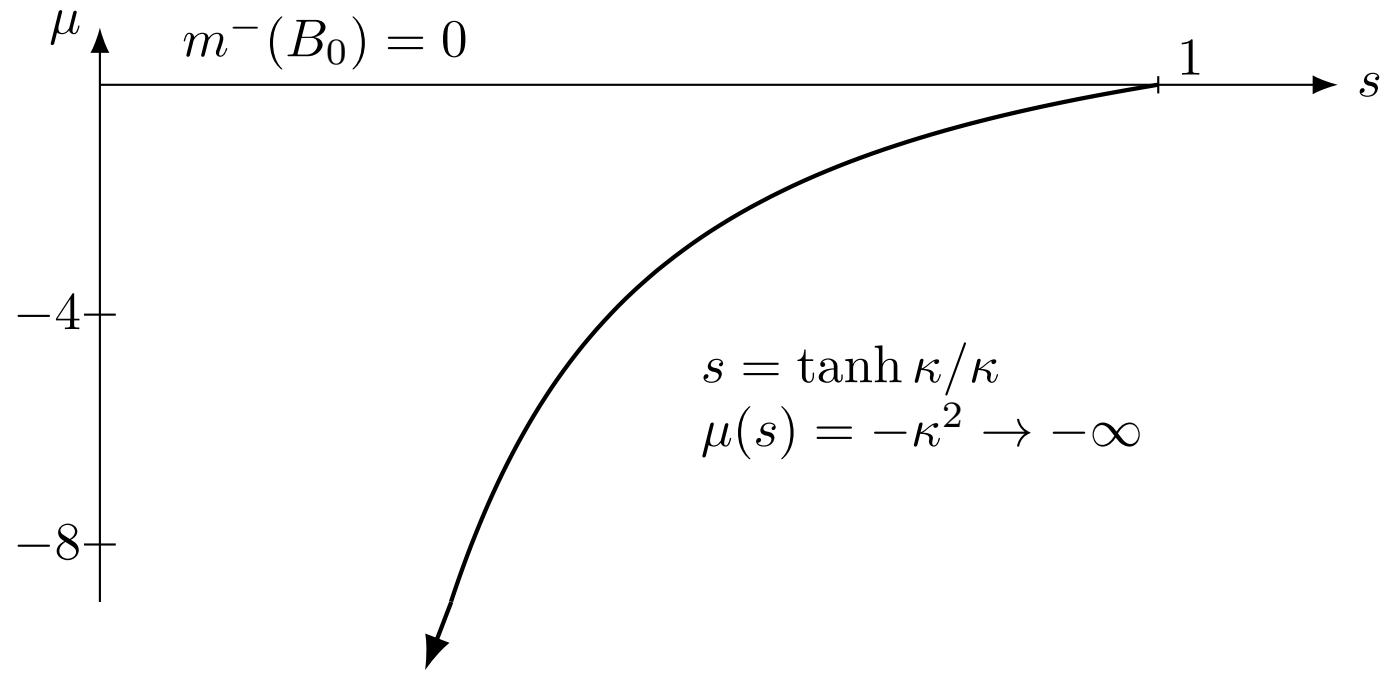
\begin{figure}[htbp]\centering
\begin{tikzpicture}[x=7cm,y=.38cm,>=Latex]
\draw[->] (0,-9)--(0,1) node[left] {$\mu$};
\draw[->] (0,0)--(1.17,0) node[right] {$s$};
\draw (1,.15)--(1,-.15) node[above right] {$1$};
\foreach \y in {-4,-8} {\draw (-.015,\y)--(.015,\y) node[left,xshift=-3pt] {$\y$};}
\draw[thick,domain=.025:3,samples=90,variable=\k]
 plot ({(exp(2*\k)-1)/(exp(2*\k)+1)/\k},{-\k*\k});
\draw[->,thick] (.332,-9)--(.307,-10.2);
\node[align=left,anchor=west] at (.55,-5.5)
 {$s=\tanh\kappa/\kappa$\\$\mu(s)=-\kappa^2\to-\infty$};
\node[anchor=west] at (.06,.75) {$m^-(B_0)=0$};
\end{tikzpicture}
\caption{The negative eigenvalue branch $\mu(s)=-\kappa(s)^2$ of $B_s$ in Proposition~\ref{prop:RobinObstruction}, which escapes to $-\infty$ as $s\downarrow0$ without crossing zero.}
\label{fig:robin}
\end{figure}

\begin{proposition}[Robin--Dirichlet topological obstruction]\label{prop:RobinObstruction}
On $\cH=L^2(0,1)$, consider the operators $B_su=-u''$ with domain
$$
 D(B_s)=\big\{u\in H^2(0,1)\ \big|\ u(1)=0,\ u(0)+su'(0)=0\big\},\qquad s\in[0,s_0],\quad 0<s_0<1.
$$
Then $s\mapsto B_s$ is a gap-continuous path of invertible self-adjoint operators with compact resolvent. Nevertheless, $m^-(B_0)=0$ while $m^-(B_s)=1$ for all $s\in(0,s_0]$, and $s\mapsto\Psi(B_s)$ is discontinuous at $s=0$.
\end{proposition}
\begin{proof}
For $z=i$, the fundamental solutions of $-u''=iu$ satisfying $u(0)+su'(0)=0$ and $u(1)=0$ vary continuously with $s\in[0,s_0]$, so that the Green kernel of $(B_s-iI)^{-1}$ converges uniformly on $[0,1]^2$, proving gap continuity. Also, any element in $\ker B_s$ is of the form $u(t)=c(1-t)$, for which $u(0)+su'(0)=c(1-s)=0$; since $s\in[0,s_0]\subset[0,1)$, we get $c=0$ and hence $\ker B_s=\{0\}$ for all $s\in[0,s_0]$.

At $s=0$ we have the Dirichlet operator $B_0>0$, so $m^-(B_0)=0$. For $s\in(0,s_0]$, the quadratic form of $B_s$ on $\{u\in H^1(0,1):u(1)=0\}$ is $q_s(u,u)=\int_0^1|u'|^2\,dt-s^{-1}|u(0)|^2$. Looking for negative eigenvalues $\mu=-\kappa^2<0$ ($\kappa>0$) gives $u(t)=c\sinh(\kappa(1-t))$, and the left boundary condition reads
\[
 s=\frac{\tanh\kappa}{\kappa}.
\]
Since $\kappa\mapsto(\tanh\kappa)/\kappa$ is strictly decreasing from $1$ to $0$ on $(0,\infty)$, for each $s\in(0,s_0]$ there is a unique negative eigenvalue $\mu(s)=-\kappa(s)^2$, which escapes to $-\infty$ as $s\downarrow0$ (see Figure~\ref{fig:robin}). Thus $m^-(B_s)=1$ for all $s\in(0,s_0]$. Choosing a normalized eigenvector $e_s$ for $\mu(s)$ and a constant $c>0$ such that $\Psi(B_0)\ge cI$, we find
\[
 \norm{\Psi(B_s)-\Psi(B_0)}
 \ge \bigl|\ip{(\Psi(B_s)-\Psi(B_0))e_s}{e_s}\bigr|
 \ge c-\Psi(\mu(s))\longrightarrow c+1>0\qquad\text{as }s\downarrow0,
\]
which shows that $s\mapsto\Psi(B_s)$ is discontinuous at $s=0$.
\end{proof}

By contrast, when the Robin coefficient $h=-s^{-1}$ varies over a compact interval, the form domain $V_R=\{u\in H^1(0,1):u(1)=0\}$ is independent of $h$. For example, if $0<\alpha<\pi/2$ and $\beta>0$, critical points of the $C_2$-invariant functional
$$
 E_h(u)=\frac{1}{2}\int_0^1(u'^2-\alpha^2u^2)\,dt+\frac{h}{2}u(0)^2+\frac{\beta}{4}\int_0^1u^4\,dt,\qquad u\in V_R,
$$
solve $-u''-\alpha^2u+\beta u^3=0$ with $u(1)=0$ and $u'(0)=hu(0)$. At $h_*=-\alpha\cot\alpha$, $\ker D^2E_{h_*}(0)$ is spanned by $v(t)=\sin(\alpha(1-t))$ and $\partial_hD^2E_h(0)|_{h_*}[v,v]=v(0)^2=\sin^2\alpha>0$; thus Corollary~\ref{cor:simple} gives a local branch $u(s)=sv+O(s^3)$ with $h(s)=h_*-\beta\big(\int_0^1v^4\,dt/v(0)^2\big)s^2+O(s^4)$ and trivial stabilizer in $C_2$.

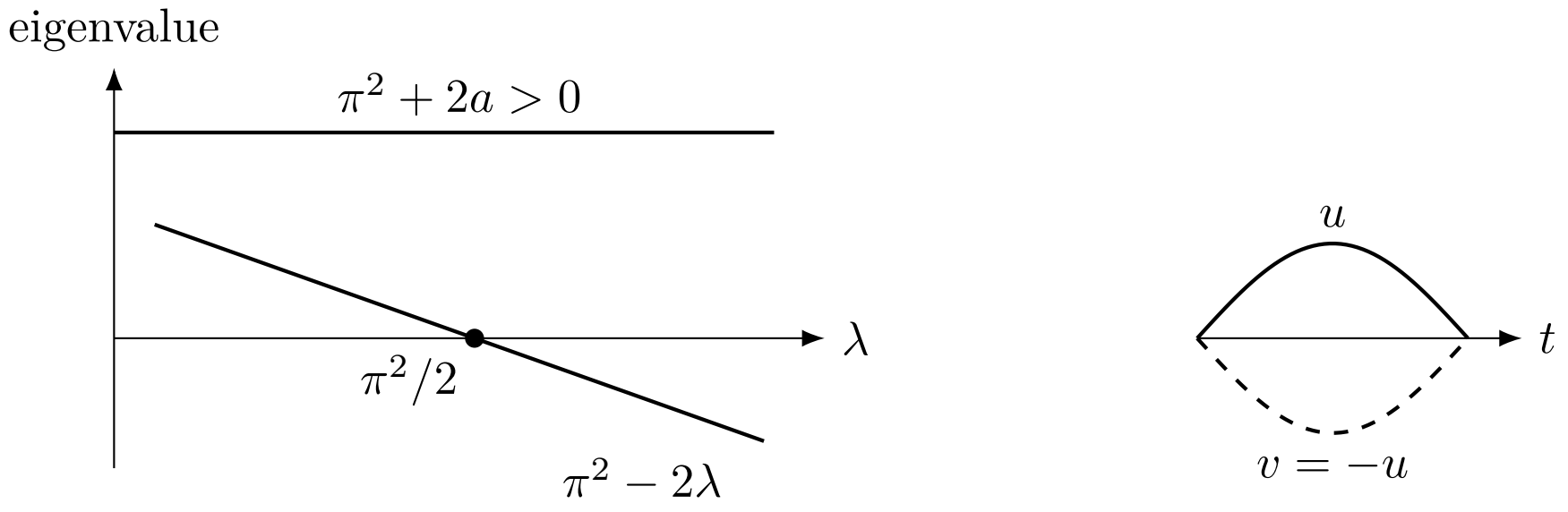
\begin{figure}[htbp]\centering
\begin{tikzpicture}[>=Latex,x=1.5cm,y=.8cm]
\begin{scope}
\draw[->] (0,0)--(3.5,0) node[right] {$\lambda$};
\draw[->] (0,-1.2)--(0,2.5) node[above] {eigenvalue};
\draw[thick] (0,1.9)--(3.25,1.9);
\node[above] at (1.7,1.9) {$\pi^2+2a>0$};
\draw[thick] (.2,1.05)--(3.2,-.95);
\fill (1.775,0) circle(2pt);
\node[below left] at (1.775,0) {$\pi^2/2$};
\node[below] at (2.6,-.95) {$\pi^2-2\lambda$};
\end{scope}
\begin{scope}[xshift=8cm,x=2cm,y=.7cm]
\draw[->] (0,0)--(1.2,0) node[right] {$t$};
\draw[thick,domain=0:1,samples=50] plot (\x,{sin(180*\x)});
\draw[thick,dashed,domain=0:1,samples=50] plot (\x,{-sin(180*\x)});
\node[above] at (.5,1) {$u$};\node[below] at (.5,-1) {$v=-u$};
\end{scope}
\end{tikzpicture}
\caption{Left: the symmetric channel $\cH^G$ remains positive definite while the antisymmetric channel crosses zero at $\lambda_1=\pi^2/2$. Right: the leading antisymmetric mode $v=-u$ bifurcating at $\lambda_1$.}
\label{fig:channels}
\end{figure}

In a similar way, if we consider $H^1_0(0,1;\R^2)$ with $G=C_2$ acting by $(u,v)\mapsto(v,u)$ and set $s=(u+v)/\sqrt{2}$, $z=(u-v)/\sqrt{2}$, then the action becomes $(s,z)\mapsto(s,-z)$. For the functional
$$
 E_\lambda(s,z)=\frac{1}{2}\int_0^1(s'^2+z'^2)\,dt+a\int_0^1s^2\,dt-\lambda\int_0^1z^2\,dt+\frac{\beta}{4}\int_0^1(s^2+z^2)^2\,dt \qquad (a\ge0,\ \beta>0),
$$
the restriction $L_+=-d^2/dt^2+2a$ to $\cH^G=\{z=0\}$ has positive eigenvalues $\{n^2\pi^2+2a\}_{n\ge1}$ for all $\lambda$, whence $\nu_G(\lambda)=0$ and $B_{\{e\}}=0$. Meanwhile, the antisymmetric block $L_-(\lambda)=-d^2/dt^2-2\lambda$ crosses zero at $\lambda_n=n^2\pi^2/2$ ($n\ge1$) with simple kernel $\sin(n\pi t)$ and local flow $-1$. Pairing the first Euler--Lagrange equation $-s''+2as+\beta(s^2+z^2)s=0$ with $s\in H^1_0(0,1)$ forces $s\equiv0$, and Theorem~\ref{thm:budget} together with Corollary~\ref{cor:simple} yields symmetry-breaking bifurcation at every $\lambda_n$ into purely antisymmetric solutions $u=-v$ ($z\ne0$) with trivial stabilizer (Figure~\ref{fig:channels}).

\section{Tangency and bifurcation of Lagrangian intersections}\label{sec:intersections}\label{subsec:lagrangian}

Theorem A also admits a natural geometric reformulation for the bifurcation of equivariant Lagrangian intersections in a Weinstein generating chart, complementing the global generating-family approach of Ciriza and Pejsachowicz \cite{CirizaPejsachowicz}.

Consider a $2d$-dimensional symplectic manifold $(M,\omega)$ endowed with a symplectic action of a finite group $G$, and let $\Lambda_{0,\lambda},\Lambda_{1,\lambda}\subset M$ be two $G$-invariant Lagrangian submanifolds meeting along a continuous curve of $G$-fixed points $z_\lambda\in\Lambda_{0,\lambda}\cap\Lambda_{1,\lambda}$, $\lambda\in I=[a,b]$. Suppose that around $z_\lambda$ we have a local $G$-equivariant Weinstein chart modeled on $T^*W$, where $W$ is a $d$-dimensional orthogonal $G$-representation with the cotangent action $g(u,p)=(gu,p\circ g^{-1})$; in other words, a continuous family of $G$-equivariant symplectomorphisms $\chi_\lambda:U_\lambda\subset M\to B_\rho(0)\times B_\rho(0)\subset T^*W$ such that
$$
 \chi_\lambda(z_\lambda)=(0,0),\qquad
 \chi_\lambda(\Lambda_{0,\lambda}\cap U_\lambda)=\big\{(u,0)\ \big|\ u\in B_\rho(0)\big\},\qquad
 \chi_\lambda(\Lambda_{1,\lambda}\cap U_\lambda)=\big\{(u,dh_\lambda(u))\ \big|\ u\in B_\rho(0)\big\},
$$
for a $G$-invariant $C^2$ family $h_\lambda:B_\rho(0)\subset W\to\R$ satisfying (H1)--(H2) and $dh_\lambda(0)=0$. If we set $S_\lambda=D^2h_\lambda(0)\in\mathcal{L}(W)$, then the tangent spaces at $(0,0)$ in the chart are $W\times\{0\}$ and $\gr S_\lambda$, whose intersection corresponds to $\ker S_\lambda$ and has $K$-fixed dimension $\nu_K(\lambda)=\dim\ker S_\lambda^K$.

\begin{corollary}[Lagrangian intersection bifurcation with prescribed symmetry]\label{cor:lagrangian}
In the setting above, let $H\le G$ be a subgroup such that $S_a^H,S_b^H\in\mathcal{L}(W^H)$ be invertible and
$$
 s_H=m^-(S_a^H)-m^-(S_b^H)=\Mas\big(W^H\times\{0\},\gr S_\lambda^H;[a,b]\big)\ne0.
$$
\begin{enumerate}[(i)]
\item There exist $t\in(a,b)$, $\lambda_j\to t$, and $H$-invariant intersection points $z_j\in(\Lambda_{0,\lambda_j}\cap\Lambda_{1,\lambda_j})^H\setminus\{z_{\lambda_j}\}$ such that $\chi_{\lambda_j}(z_j)\to(0,0)$.
\item If $\Sigma_H=\{t\in(a,b):\dim\ker S_t^H>0\}$ is finite and $|s_H|>B_H$, then at least $\lceil(|s_H|-B_H)/M_H\rceil$ distinct instants $t\in\Sigma_H\setminus\cE_H$ admit such sequences of intersection points with $\Stab_G(z_j)=H$, where $\cE_H,B_H,M_H$ are defined from $S_\lambda$ as in \eqref{eq:exceptional}--\eqref{eq:BH_MH}.
\end{enumerate}
\end{corollary}
\begin{proof}
In the coordinates given by $\chi_\lambda$, points $(u,0)\in\chi_\lambda(\Lambda_{0,\lambda})\cap\chi_\lambda(\Lambda_{1,\lambda})$ are in one-to-one correspondence with critical points $dh_\lambda(u)=0$ in $B_\rho(0)\subset W$. The statement then follows by applying Theorem A to $h_\lambda$ on $W$ and taking $z_j=\chi_{\lambda_j}^{-1}(u_j,0)$.
\end{proof}

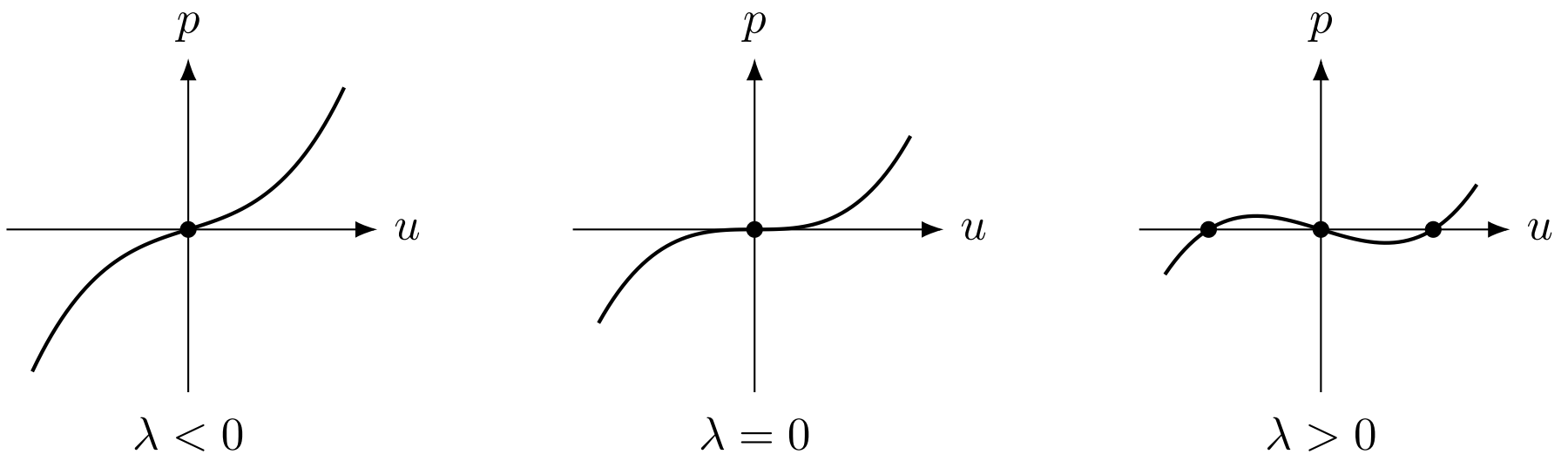
\begin{figure}[htbp]\centering
\begin{tikzpicture}[>=Latex,x=1.15cm,y=.65cm]
\foreach \shift/\lam/\txt in {0/-.55/{\lambda<0},4.3/0/{\lambda=0},8.6/.55/{\lambda>0}} {
\begin{scope}[xshift=\shift cm]
\draw[->] (-1.2,0)--(1.25,0) node[right] {$u$};
\draw[->] (0,-1.9)--(0,2) node[above] {$p$};
\draw[thick,domain=-1.03:1.03,samples=65]
 plot (\x,{\x*\x*\x-\lam*\x});
\fill (0,0) circle(1.8pt);
\node[below] at (0,-2) {$\txt$};
\end{scope}}
\begin{scope}[xshift=8.6cm]
\fill ({sqrt(.55)},0) circle(1.8pt);
\fill ({-sqrt(.55)},0) circle(1.8pt);
\end{scope}
\end{tikzpicture}
\caption{Lagrangian intersection bifurcation for the graph $p=u^3-\lambda u$ and the zero section $p=0$: the non-zero Maslov index change across the tangency at $\lambda=0$ produces a symmetric pair of intersection points with trivial $C_2$-stabilizer for $\lambda>0$.}
\label{fig:tangency}
\end{figure}

Here we follow the relative Maslov index convention of Robbin and Salamon \cite{RobbinSalamon} in which the first Lagrangian subspace $W^H\times\{0\}$ is held fixed and the moving graph $\gr S_\lambda^H$ has crossing form $\ip{\dot S_\lambda^H v}{v}$ on $\ker S_\lambda^H$, so that $\Mas(W^H\times\{0\},\gr S_\lambda^H;[a,b])=\sflo(S^H;[a,b])=s_H$. Observe that tangency ($\ker S_t^H\ne\{0\}$) without a Maslov index jump does not force bifurcation: on $T^*\R$ with $C_2$ acting by $-1$, $h_\lambda(u)=\frac{1}{4}u^4-\frac{1}{2}\lambda u^2$ has $s_{\{e\}}=-1$ and bifurcates into $(\pm\sqrt{\lambda},0)$ for $\lambda>0$ (Figure~\ref{fig:tangency}), whereas $\widetilde{h}_\lambda(u)=\frac{1}{4}u^4+\frac{1}{2}\lambda^2u^2$ becomes tangent at $\lambda=0$ with $s_{\{e\}}=0$ and admits no non-trivial real intersections.

\section{Central configurations and periodic motions in the \texorpdfstring{$N$}{N}-body problem}\label{sec:mechanics}

We now turn to the Newtonian $N$-body problem. In Subsection~\ref{subsec:cc} we prove a general computational lemma showing that the restricted shape Morse indices and nullities of a planar central configuration can be computed directly from a rank-one rotational lift on the unconstrained configuration space $\R^{2N}$. After illustrating this criterion on Palmore's centered equilateral triangle in Subsection~\ref{subsec:triangle}, we establish Theorem~B for the nine-body interlaced-square family in Subsection~\ref{subsec:ninebody}, and conclude in Subsection~\ref{subsec:periodic} with a bifurcation result for spatial relative-periodic orbits of a vertically oscillating ring.

\subsection{The constrained shape Hessian and the rotational lift}\label{subsec:cc}

Let $m_0(\lambda),\ldots,m_{N-1}(\lambda)>0$ be smooth point masses in the plane $\R^2$, $\lambda\in[a,b]$. On the collision-free configuration space in $\R^{2N}$, we denote the Newtonian potential, the moment of inertia, and the mass matrix respectively by
$$
 U_\lambda(q)=\sum_{0\le i<j\le N-1}\frac{m_i(\lambda)m_j(\lambda)}{|q_i-q_j|},\qquad
 I_\lambda(q)=\sum_{i=0}^{N-1}m_i(\lambda)|q_i|^2,\qquad
 M_\lambda=\diag\big(m_0(\lambda)I_2,\dots,m_{N-1}(\lambda)I_2\big).
$$
Planar central configurations are critical points of the restriction of $U_\lambda$ to the inertia ellipsoid
\[
 \mathcal{S}_\lambda=\Bigl\{q\in\R^{2N}:\sum_{i=0}^{N-1}m_i(\lambda)q_i=0,\ I_\lambda(q)=1\Bigr\},
\]
considered modulo the planar rotation action of $SO(2)$ generated by $J=\diag(J_2,\dots,J_2)$, where $J_2=\bigl(\begin{smallmatrix}0&-1\\1&0\end{smallmatrix}\bigr)$ is the counterclockwise quarter-turn on $\R^2$. At an unnormalized central configuration $q_\lambda$ with center of mass at the origin, the central-configuration equation reads
\[
 \nabla U_\lambda(q_\lambda)=-\omega_\lambda^2M_\lambda q_\lambda,\qquad
 \omega_\lambda^2=\frac{U_\lambda(q_\lambda)}{I_\lambda(q_\lambda)}>0,
\]
where the expression for $\omega_\lambda^2$ follows from Euler's homogeneous function relation $\ip{\nabla U_\lambda(q)}{q}=-U_\lambda(q)$.

Suppose now that $q_\lambda$ ($\lambda\in[a,b]$) be a smooth curve of non-collinear planar central configurations fixed by a finite group $G$ ($gq_\lambda=q_\lambda$) which acts orthogonally on $\R^{2N}$ through a homomorphism $g\mapsto(\rho(g),\sigma_g)\in O(2)\times S_N$ preserving the masses ($m_{\sigma_g(i)}(\lambda)=m_i(\lambda)$), i.e.
$$
 (gq)_i=\rho(g)q_{\sigma_g^{-1}(i)}.
$$
Since $q_\lambda$ is not collinear, the $SO(2)$-action at $\bar q_\lambda=q_\lambda/\sqrt{I_\lambda(q_\lambda)}\in\mathcal S_\lambda$ is free. Cutting $\mathcal S_\lambda$ orthogonally to the rotational tangent vector $J\bar q_\lambda$ with respect to the mass inner product $\ip{\xi}{\eta}_{M_\lambda}=\xi^TM_\lambda\eta$ defines a $(2N-4)$-dimensional local $G$-invariant shape slice
\[
 \Sigma_\lambda=\bigl\{q\in\mathcal S_\lambda:\ip{M_\lambda J\bar q_\lambda}{q}=0,\ \norm{q-\bar q_\lambda}<\delta_\lambda\bigr\},
\]
whose tangent space at $\bar q_\lambda$ is the $M_\lambda$-orthogonal complement
\[
 W_\lambda=T_{\bar q_\lambda}\Sigma_\lambda=\bigl(\mathcal{T}_{\mathrm{tr}}\oplus\R q_\lambda\oplus\R Jq_\lambda\bigr)^{\perp_{M_\lambda}},
\]
where $\mathcal{T}_{\mathrm{tr}}=\{(\xi_0,\dots,\xi_{N-1}):\xi_0=\dots=\xi_{N-1}\in\R^2\}$ is the two-dimensional space of uniform translations. The radial projection chart
\begin{equation}\label{eq:radial_chart}
 \phi_\lambda:B_{r_\lambda}(0)\subset W_\lambda\longrightarrow\Sigma_\lambda,\qquad
 \phi_\lambda(\eta)=\frac{\bar q_\lambda+\eta}{\sqrt{I_\lambda(\bar q_\lambda+\eta)}},
\end{equation}
is $G$-equivariant and satisfies $d\phi_\lambda(0)=\mathrm{id}_{W_\lambda}$. By $(-1)$-homogeneity of $U_\lambda$, the Hessian of the pulled-back potential $U_\lambda\circ\phi_\lambda$ at $\eta=0\in W_\lambda$ is a positive scalar multiple of the constrained Hessian $K_\lambda|_{W_\lambda}$:
\begin{equation}\label{eq:K_def}
 D^2(U_\lambda\circ\phi_\lambda)(0)[\xi,\zeta]=I_\lambda(q_\lambda)^{3/2}\ip{K_\lambda\xi}{\zeta}\quad(\xi,\zeta\in W_\lambda),\qquad
 K_\lambda=D^2U_\lambda(q_\lambda)+\omega_\lambda^2M_\lambda.
\end{equation}
Furthermore, since $M_\lambda$ and the four-dimensional symmetry subspace $\mathcal{T}_{\mathrm{tr}}\oplus\R q_\lambda\oplus\R Jq_\lambda$ vary smoothly with $\lambda\in[a,b]$, the $M_\lambda$-orthogonal projections onto $W_\lambda$ induce smooth $G$-equivariant Kato isomorphisms $\Upsilon_\lambda:W_a\xrightarrow{\sim}W_\lambda$ (constructed as in Section~\ref{subsec:regularforms}), placing the family $f_\lambda=U_\lambda\circ\phi_\lambda\circ\Upsilon_\lambda$ on the fixed Hilbert space $W_a$ without altering the restricted Morse indices or nullities.

To avoid constructing explicit parameter-dependent bases for the $(2N-4)$-dimensional slice tangent space $W_\lambda$, we introduce the rank-one rotational lift on the full configuration space $\R^{2N}$:
\begin{equation}\label{eq:liftedmechanical}
 \mathscr{Q}_\lambda=D^2U_\lambda(q_\lambda)+\omega_\lambda^2M_\lambda+\frac{\omega_\lambda^2}{I_\lambda(q_\lambda)}u_\lambda u_\lambda^T,\qquad
 u_\lambda=M_\lambda Jq_\lambda.
\end{equation}
Notice that $\rho(g)J_2=\det(\rho(g))J_2\rho(g)=\pm J_2\rho(g)$ for every $g\in G$, whence $gu_\lambda=\pm u_\lambda$, which shows that both $\R Jq_\lambda$ and the rank-one matrix $u_\lambda u_\lambda^T$ are $G$-equivariant. We also set $\widetilde{K}_\lambda=\omega_\lambda^{-2}\mathscr{Q}_\lambda$. For any vector $\xi=(\xi_0,\dots,\xi_{N-1})\in\R^{2N}$, a straightforward computation gives
\begin{equation}\label{eq:pairhessian}
 D^2U_\lambda(q)[\xi,\xi]=\sum_{0\le i<j\le N-1}\frac{m_im_j}{r_{ij}^3}
 \left(3\frac{\ip{q_i-q_j}{\xi_i-\xi_j}^2}{r_{ij}^2}-|\xi_i-\xi_j|^2\right),\qquad r_{ij}=|q_i-q_j|.
\end{equation}

\begin{lemma}[Computing restricted shape indices via the rotational lift]\label{lem:mechanicalinertia}
For every subgroup $H\le G$, the Morse index and nullity of the constrained shape Hessian $K_\lambda|_{W_\lambda^H}$ coincide with the Morse index and nullity of the lifted matrix $\mathscr{Q}_\lambda$ (and of $\widetilde{K}_\lambda=\omega_\lambda^{-2}\mathscr{Q}_\lambda$) restricted to the full configuration fixed space $(\R^{2N})^H$:
\[
 m^-\bigl(K_\lambda|_{W_\lambda^H}\bigr)=m^-\bigl(\mathscr{Q}_\lambda|_{(\R^{2N})^H}\bigr),\qquad
 \dim\ker\bigl(K_\lambda|_{W_\lambda^H}\bigr)=\dim\ker\bigl(\mathscr{Q}_\lambda|_{(\R^{2N})^H}\bigr).
\]
\end{lemma}
\begin{proof}
With respect to the mass scalar product $\ip{\xi}{\eta}_{M_\lambda}=\xi^TM_\lambda\eta$, we have the $M_\lambda$-orthogonal, $G$-invariant decomposition
$$
 \R^{2N}=\mathcal{T}_{\mathrm{tr}}\oplus\R q_\lambda\oplus\R Jq_\lambda\oplus W_\lambda.
$$
By translation invariance, $(-1)$-homogeneity, and rotational invariance of $U_\lambda$, we readily compute $D^2U_\lambda(q_\lambda)|_{\mathcal{T}_{\mathrm{tr}}}=0$, $D^2U_\lambda(q_\lambda)q_\lambda=-2\nabla U_\lambda(q_\lambda)=2\omega_\lambda^2M_\lambda q_\lambda$, and $D^2U_\lambda(q_\lambda)(Jq_\lambda)=J\nabla U_\lambda(q_\lambda)=-\omega_\lambda^2M_\lambda Jq_\lambda$. Consequently, on the first three summands the operator $M_\lambda^{-1}K_\lambda=\omega_\lambda^2I+M_\lambda^{-1}D^2U_\lambda(q_\lambda)$ acts by
$$
 M_\lambda^{-1}K_\lambda\big|_{\mathcal{T}_{\mathrm{tr}}}=\omega_\lambda^2I_2>0,\qquad
 M_\lambda^{-1}K_\lambda(q_\lambda)=3\omega_\lambda^2q_\lambda>0,\qquad
 M_\lambda^{-1}K_\lambda(Jq_\lambda)=0.
$$
The rank-one term $\frac{\omega_\lambda^2}{I_\lambda(q_\lambda)}u_\lambda u_\lambda^T$ vanishes on $\mathcal{T}_{\mathrm{tr}}\oplus\R q_\lambda\oplus W_\lambda$ and moves the eigenvalue on $\R Jq_\lambda$ from $0$ to $\omega_\lambda^2>0$. Thus $\mathscr{Q}_\lambda$ is positive definite on $\mathcal{T}_{\mathrm{tr}}\oplus\R q_\lambda\oplus\R Jq_\lambda$ and agrees with $K_\lambda$ on $W_\lambda$; restricting to $(\R^{2N})^H$ and applying Sylvester's law of inertia finishes the proof.
\end{proof}

From Theorem A and Lemma~\ref{lem:mechanicalinertia} we immediately get the following criterion for central configurations.

\begin{corollary}[Fixed-space test for central configurations]\label{cor:cc}
Let $q_\lambda$, $\lambda\in[a,b]$, be a smooth curve of collision-free non-collinear planar central configurations with finite isotropy group $G$, and let $H\le G$ be a subgroup such that $\mathscr{Q}_a|_{(\R^{2N})^H}$ and $\mathscr{Q}_b|_{(\R^{2N})^H}$ be invertible with
$$
 s_H=m^-\big(\mathscr{Q}_a|_{(\R^{2N})^H}\big)-m^-\big(\mathscr{Q}_b|_{(\R^{2N})^H}\big)\ne0.
$$
Then non-trivial $H$-invariant central configurations bifurcate from $q_\lambda$ (up to rotations and scaling) at some instant in $(a,b)$; moreover, if $\Sigma_H$ is finite and $|s_H|>B_H$, at least $\lceil(|s_H|-B_H)/M_H\rceil$ distinct instants in $(a,b)$ admit bifurcating central configurations with exact stabilizer $H$.
\end{corollary}

\subsection{Benchmark example: central displacement from an equilateral triangle}\label{subsec:triangle}

As a quick warm-up before the nine-body problem, let us see how Corollary~\ref{cor:cc} and Corollary~\ref{cor:simple} retrieve Meyer and Schmidt's isosceles bifurcation \cite{MeyerSchmidt1988} at Palmore's critical mass \cite{Palmore1976} (cf.\ also \cite{SantosLeandro2026}).

\begin{figure}[htbp]\centering
\begin{tikzpicture}[>=Latex,scale=1.25]
\begin{scope}[xshift=-2.3cm]
\draw (0:1)--(120:1)--(240:1)--cycle;
\foreach \t in {0,120,240} {\fill (\t:1) circle(1.7pt);}
\fill (0,0) circle(2.2pt);
\draw[dashed] (-1.1,0)--(1.2,0);
\node[below] at (0,-1.15) {centered equilateral ($D_3$)};
\node[above] at (0,0) {$\mu$};
\end{scope}
\draw[->] (-.8,0)--(.6,0) node[midway,above] {$\mu=\mu_*$};
\begin{scope}[xshift=2cm]
\draw (.95,0)--(-.5,.9)--(-.5,-.9)--cycle;
\fill (.95,0) circle(1.7pt);\fill (-.5,.9) circle(1.7pt);
\fill (-.5,-.9) circle(1.7pt);\fill (.12,0) circle(2.2pt);
\draw[dashed] (-1,0)--(1.2,0);
\draw[->] (.12,.17)--(.52,.17);
\node[below] at (0,-1.15) {isosceles branch ($H\cong C_2$)};
\end{scope}
\end{tikzpicture}
\caption{Symmetry breaking from $D_3$ to $H\cong C_2$ at Palmore's critical mass $\mu_*=(81+64\sqrt{3})/249$, accompanied by displacement of the central body along the reflection axis.}
\label{fig:triangle}
\end{figure}
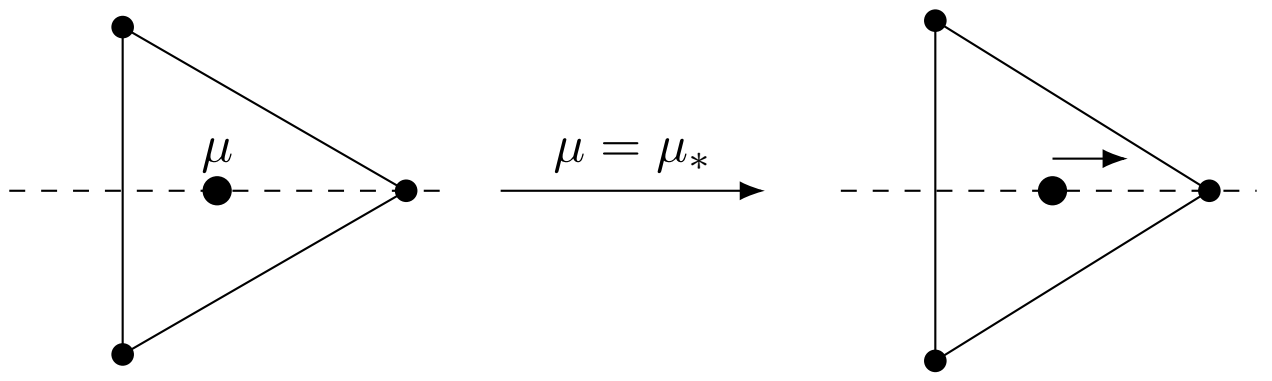

Consider four bodies with masses $(m_0,m_1,m_2,m_3)=(\mu,1,1,1)$ at the centered equilateral triangle positions
\[
 q_0=(0,0),\qquad q_1=(1,0),\qquad q_2=\bigl(-\tfrac{1}{2},\tfrac{\sqrt{3}}{2}\bigr),\qquad q_3=\bigl(-\tfrac{1}{2},-\tfrac{\sqrt{3}}{2}\bigr),
\]
for which $I=3$ and $\omega^2=\mu+1/\sqrt{3}$. The branch isotropy group is $D_3$, and we take $H\cong C_2$ to be reflection across the horizontal axis coupled with swapping bodies $2$ and $3$. The $H$-fixed shape slice $W^H$ is two-dimensional and is parametrized by $(x,y)\in\R^2$ via
\[
 \xi_0=\bigl((-3x+\sqrt{3}y)/\mu,0\bigr),\qquad
 \xi_1=(x-\sqrt{3}y,0),\qquad
 \xi_2=(x,y),\qquad
 \xi_3=(x,-y),
\]
which satisfy $\sum_{i=0}^3m_i\xi_i=0$, $\ip{Mq}{\xi}=0$, and $\ip{MJq}{\xi}=0$. By a direct computation using \eqref{eq:pairhessian}, the quadratic form $C_\mu[\xi,\xi]=(D^2U(q)+\omega^2M)[\xi,\xi]$ on $W^H$ is represented by a $2\times2$ matrix $C_\mu$ satisfying
$$
 (C_\mu)_{11}=\frac{9}{2}\mu+\sqrt{3}+18+\frac{3\sqrt{3}+27/2}{\mu}>0,\qquad
 \det C_\mu=(15\sqrt{3}-54)\frac{\mu+3}{\mu}(\mu-\mu_*),
$$
where $\mu_*=(81+64\sqrt{3})/249\approx0.77048695$. Since $15\sqrt{3}-54<0$, the index $m^-(C_\mu)$ changes from $0$ for $\mu<\mu_*$ to $1$ for $\mu>\mu_*$, while $W^{D_3}=\{0\}$ gives $\nu_{D_3}(\mu)=0$ and $B_H=0$. Moreover, setting $\xi_0=0$ forces $y=\sqrt{3}x$, where $(1,\sqrt{3})C_\mu(1,\sqrt{3})^T=18\mu+6\sqrt{3}>0$. Consequently, Corollary~\ref{cor:cc} and Corollary~\ref{cor:simple} yield a smooth local branch of non-equilateral isosceles central configurations bifurcating at $\mu_*$ with exact stabilizer $H$ and first-order displacement of the central body along the reflection axis (Figure~\ref{fig:triangle}).

\subsection{The nine-body interlaced-square family and proof of Theorem B}\label{subsec:ninebody}

We now turn to our main mechanical result: simultaneous symmetry breaking and central-body displacement in the Newtonian nine-body system formed by two interlaced squares of satellites and a central mass.

\subsubsection{Exact central-configuration branch and \texorpdfstring{$D_4$}{D4} symmetry}\label{subsubsec:nine_setup}

For a radius parameter $a$ near $a_0=177/250=0.708$, consider nine planar bodies at the positions
\begin{equation}\label{eq:ninepositions}
 \begin{split}
 q_0&=(0,0),\qquad
 (q_1,q_2,q_3,q_4)=\bigl((1,0),(0,1),(-1,0),(0,-1)\bigr),\\
 (q_5,q_6,q_7,q_8)&=\bigl((a,a),(-a,a),(-a,-a),(a,-a)\bigr).
 \end{split}
\end{equation}
Here bodies $1,2,3,4$ form the first square of radius $1$ on the coordinate axes, and bodies $5,6,7,8$ form the second square of radius $R(a)=\sqrt{2}a$ on the diagonals (since $a_0=0.708>1/\sqrt{2}\approx0.70710678$, we have $R(a)>1$, so the two squares are interlaced). We introduce the auxiliary geometric functions
\begin{equation}\label{eq:ninegeom}
 \begin{split}
 R(a)&=\sqrt{2}a,\qquad d_\pm(a)=\sqrt{1\pm2a+2a^2},\qquad S=\frac{1+2\sqrt{2}}{4},\\
 F(a)&=\frac{2(1-a)}{d_-(a)^3}+\frac{2(1+a)}{d_+(a)^3},\qquad
 G(a)=\frac{2-1/a}{d_-(a)^3}+\frac{2+1/a}{d_+(a)^3},
 \end{split}
\end{equation}
and fix the second-square satellite mass $\beta_0>0$ once and for all by the exact algebraic value at $a_0=177/250$:
\begin{equation}\label{eq:ninebeta}
 \beta_0=\frac{G(a_0)-S-150\bigl(1-R(a_0)^{-3}\bigr)}{F(a_0)-SR(a_0)^{-3}}
 =0.3650570181871535\ldots
\end{equation}
Holding all eight satellite masses fixed at $(m_1,\dots,m_8)=(1,1,1,1,\beta_0,\beta_0,\beta_0,\beta_0)$, we define the central mass $\mu(a)$ along the branch by
\begin{equation}\label{eq:ninemasses}
 \mu(a)=\frac{G(a)-S-\beta_0\bigl(F(a)-SR(a)^{-3}\bigr)}{1-R(a)^{-3}},
 \qquad (m_0,m_1,\dots,m_8)=(\mu(a),1,1,1,1,\beta_0,\beta_0,\beta_0,\beta_0).
\end{equation}
By reflection symmetry across the coordinate axes and diagonals, all tangential forces on the satellites and the net force on the central body $q_0=(0,0)$ vanish identically. The radial acceleration multipliers $-\nabla_{q_i}U(q)/(m_iq_i)$ on the first and second squares are $\mu(a)+S+\beta_0F(a)$ and $\mu(a)R(a)^{-3}+\beta_0SR(a)^{-3}+G(a)$, which coincide for every $a$ by \eqref{eq:ninemasses}. Consequently, \eqref{eq:ninepositions}--\eqref{eq:ninemasses} defines a real-analytic family of central configurations with
\begin{equation}\label{eq:nineomega}
 \omega^2(a)=\mu(a)+S+\beta_0F(a)=\mu(a)R(a)^{-3}+\beta_0SR(a)^{-3}+G(a),\qquad
 I(a)=4+8a^2\beta_0,
\end{equation}
satisfying $\mu(a_0)=150$ at $a_0=177/250$. As proved in Lemma~\ref{lem:ninecertificate} below, $\mu'(a)<0$ throughout the interval $\mathcal{J}=[177/250,14161/20000]$, so the central mass $\mu\in[\mu(a_1),150]$ serves as a smooth, strictly monotone parameter along the branch.

Let $G=D_4=\langle r,s\rangle$ be the dihedral group of order eight generated by the counterclockwise quarter-turn $r=J_2$ (accompanied by the cyclic vertex permutations $(1\,2\,3\,4)(5\,6\,7\,8)$) and the horizontal reflection $s=\diag(1,-1)$ (accompanied by the vertex transpositions $(2\,4)(5\,8)(6\,7)$). Up to conjugacy in $D_4$, the two reflection subgroups of order two are
$$
 H_{\mathrm{ax}}=\langle s\rangle \quad \textrm{ and } \quad H_{\mathrm{diag}}=\langle rs\rangle,
$$
which fix the horizontal axis resp.\ the main diagonal $\{x=y\}$. Notice that the linear action of $D_4$ on $\R^{18}$ does not depend on $a$ or $\mu$.

\subsubsection{Schur complement scalarization and simultaneous bifurcation mechanism}\label{subsubsec:schur_mechanism}

Before proving the rational interval bounds, let us explain the representation-theoretic reason why the two subspaces $W_a^{H_{\mathrm{ax}}}$ and $W_a^{H_{\mathrm{diag}}}$ must bifurcate at the same critical mass.

We write any variation $\xi\in\R^{18}$ as $\xi=(\xi_0,\xi_{\mathrm{sat}})\in\R^2\oplus\R^{16}$, where $\xi_0\in\R^2$ is the displacement of the central body and $\xi_{\mathrm{sat}}=(\xi_1,\dots,\xi_8)\in\R^{16}$ is the displacement vector of the eight satellites. With respect to this splitting, the lifted matrix $\mathscr{Q}(a)=\omega^2(a)\widetilde{K}(a)$ defined in \eqref{eq:liftedmechanical} reads
\begin{equation}\label{eq:Q_partition}
 \mathscr{Q}(a)=\begin{pmatrix}A_0(a)&B(a)^T\\B(a)&C(a)\end{pmatrix},
\end{equation}
with $A_0(a)=\alpha(a)I_2\in\R^{2\times 2}$, $B(a)\in\R^{16\times 2}$, and $C(a)=\mathscr{Q}(a)|_{\{\xi_0=0\}}\in\R^{16\times 16}$. Observe that for $H=H_{\mathrm{ax}},H_{\mathrm{diag}}$ we have $\dim(\R^{18})^H=9$, and hence quotienting out the $1$-dimensional $H$-fixed translation line and the radial direction gives $\dim W_a^{H_{\mathrm{ax}}}=\dim W_a^{H_{\mathrm{diag}}}=7$ (while $Jq(a)$ is anti-invariant under both reflections).

\begin{lemma}[$D_4$ scalarization of the central Schur complement and simultaneous crossing]\label{lem:simultaneous}
Assume that the satellite block $C(a)$ in \eqref{eq:Q_partition} be positive definite for all $a\in\mathcal{J}=[a_0,a_1]$, and set
\begin{equation}\label{eq:schur_def}
 S(a):=A_0(a)-B(a)^TC(a)^{-1}B(a)\in\R^{2\times2}.
\end{equation}
Then the following statements hold:
\begin{enumerate}[(i)]
\item For each $a\in\mathcal{J}$, $S(a)$ is a scalar multiple of $I_2$, namely
\begin{equation}\label{eq:schur_scalar}
 S(a)=\sigma(a)I_2,
\end{equation}
where $\sigma$ is a real-analytic function on a neighborhood of $\mathcal{J}$.
\item For every $a\in\mathcal{J}$, the full shape Morse index and nullity on $W_a$, as well as the restricted shape Morse indices and nullities on both reflection-fixed spaces $W_a^{H_{\mathrm{ax}}}$ and $W_a^{H_{\mathrm{diag}}}$, are completely determined by the sign of $\sigma(a)$:
\begin{align}
 m^-\bigl(K(a)|_{W_a}\bigr)&=2\cdot\mathbf{1}_{(-\infty,0)}(\sigma(a)),&
 \dim\ker\bigl(K(a)|_{W_a}\bigr)&=2\cdot\mathbf{1}_{\{0\}}(\sigma(a)),\label{eq:full_shape_sigma}\\
 m^-\bigl(K(a)|_{W_a^H}\bigr)&=\mathbf{1}_{(-\infty,0)}(\sigma(a)),&
 \dim\ker\bigl(K(a)|_{W_a^H}\bigr)&=\mathbf{1}_{\{0\}}(\sigma(a))\quad(H=H_{\mathrm{ax}},H_{\mathrm{diag}}).\label{eq:refl_shape_sigma}
\end{align}
Moreover, for every subgroup $K\le D_4$ strictly containing $H_{\mathrm{ax}}$ or $H_{\mathrm{diag}}$, one has $\nu_K(a)=\dim\ker(K(a)|_{W_a^K})=0$ for all $a\in\mathcal{J}$.
\item If $\sigma(a_0)<0<\sigma(a_1)$, then the singular set $\Sigma=\{a\in(a_0,a_1):\sigma(a)=0\}$ is non-empty and finite, and there exists at least one common critical parameter $a_*\in\Sigma$ (with critical central mass $\mu_*=\mu(a_*)$) such that the local spectral flows of both restricted Hessians at $a_*$ satisfy
\[
 s_{H_{\mathrm{ax}}}(a_*)=s_{H_{\mathrm{diag}}}(a_*)\ne0.
\]
In particular, central configurations with exact stabilizer $H_{\mathrm{ax}}$ resp.\ $H_{\mathrm{diag}}$ bifurcate simultaneously from the concentric branch at the same critical mass $\mu_*=\mu(a_*)$.
\end{enumerate}
\end{lemma}
\begin{proof}
\emph{Part (i):} Both $\R^2\times\{0\}$ and $\{0\}\times\R^{16}$ are invariant under the action of $D_4$ on $\R^{18}$, and on $\R^2$ the generators $r=J_2$ and $s=\diag(1,-1)$ give the irreducible $2$-dimensional representation $E$. Since $\mathscr{Q}(a)$ is symmetric and commutes with $D_4$, and since $C(a)>0$ is invertible, the $2\times2$ symmetric matrix $S(a)=A_0(a)-B(a)^TC(a)^{-1}B(a)$ commutes with both $r$ and $s$, which immediately implies $S(a)=\sigma(a)I_2$. Real analyticity of $\sigma$ on a neighborhood of $\mathcal J$ follows from the fact that $d_\pm(a)>0$, $1-R(a)^{-3}\ne0$, and $\det C(a)>0$ on $\mathcal J$.

\emph{Part (ii):} Since $C(a)>0$, for any $\xi=(\xi_0,\xi_{\mathrm{sat}})\in\R^2\oplus\R^{16}$ we can complete the square to write
\begin{equation}\label{eq:congruence}
 \ip{\mathscr{Q}(a)\xi}{\xi}=\ip{S(a)\xi_0}{\xi_0}+\ip{C(a)\big(\xi_{\mathrm{sat}}+C(a)^{-1}B(a)\xi_0\big)}{\xi_{\mathrm{sat}}+C(a)^{-1}B(a)\xi_0}.
\end{equation}
By Sylvester's law of inertia, we get $m^-(\mathscr{Q}(a))=m^-(S(a))=m^-(\sigma(a)I_2)$ and $\dim\ker\mathscr{Q}(a)=\dim\ker(\sigma(a)I_2)$, so that Lemma~\ref{lem:mechanicalinertia} yields \eqref{eq:full_shape_sigma}. Furthermore, for $H=H_{\mathrm{ax}}$ resp.\ $H_{\mathrm{diag}}$, the linear change of variables $(\xi_0,\xi_{\mathrm{sat}})\mapsto(\xi_0,\xi_{\mathrm{sat}}+C(a)^{-1}B(a)\xi_0)$ is $D_4$-equivariant and hence maps $(\R^{18})^H=E^H\oplus(\R^{16})^H$ onto itself. Because $E^{H_{\mathrm{ax}}}=\R(1,0)$ and $E^{H_{\mathrm{diag}}}=\R(1,1)$ are $1$-dimensional lines on which $S(a)=\sigma(a)I_2$ is simply multiplication by $\sigma(a)$, and $C(a)|_{(\R^{16})^H}>0$, Sylvester's law of inertia on $(\R^{18})^H$ and Lemma~\ref{lem:mechanicalinertia} give \eqref{eq:refl_shape_sigma}. Finally, any subgroup $K\le D_4$ strictly larger than $H_{\mathrm{ax}}$ or $H_{\mathrm{diag}}$ has order $4$ or $8$ and therefore contains the center $Z(D_4)=\{e,r^2\}$. Since $r^2|_E=-I_2$, we have $E^K=\{0\}$ and $(\R^{18})^K\subset\{0\}\times\R^{16}$, on which $\mathscr{Q}(a)|_{(\R^{18})^K}=C(a)|_{(\R^{16})^K}>0$; thus $\nu_K(a)=0$, $\forall\,a\in\mathcal{J}$.

\emph{Part~(iii):} Because $\sigma(a)$ is real analytic on $\mathcal{J}=[a_0,a_1]$ and $\sigma(a_0)<0<\sigma(a_1)$, $\sigma$ is not identically zero, so its zero set $\Sigma=\{a\in(a_0,a_1):\sigma(a)=0\}$ is non-empty and finite. By \eqref{eq:refl_shape_sigma}, both restricted paths $L^{H_{\mathrm{ax}}}$ and $L^{H_{\mathrm{diag}}}$ share the same finite singular set $\Sigma_{H_{\mathrm{ax}}}=\Sigma_{H_{\mathrm{diag}}}=\Sigma$. Moreover, on any isolating interval $[a_*-\delta,a_*+\delta]\subset(a_0,a_1)$ around a root $a_*\in\Sigma$, the endpoint Morse indices of $L^{H_{\mathrm{ax}}}$ and $L^{H_{\mathrm{diag}}}$ both equal $\mathbf{1}_{(-\infty,0)}(\sigma(a_*\pm\delta))$, so their local spectral flows coincide:
\[
 s_{H_{\mathrm{ax}}}(a_*)=\mathbf{1}_{(-\infty,0)}(\sigma(a_*-\delta))-\mathbf{1}_{(-\infty,0)}(\sigma(a_*+\delta))=s_{H_{\mathrm{diag}}}(a_*).
\]
As $\sum_{a_*\in\Sigma}s_{H_{\mathrm{ax}}}(a_*)=1-0=1\ne0$, there must be at least one zero $a_*\in\Sigma$ with $s_{H_{\mathrm{ax}}}(a_*)=s_{H_{\mathrm{diag}}}(a_*)\ne0$. Since $\cE_{H_{\mathrm{ax}}}=\cE_{H_{\mathrm{diag}}}=\varnothing$ by item (ii), Theorem~\ref{thm:budget}(i) applies at $\mu_*=\mu(a_*)$ to both subgroups at once.
\end{proof}

\subsubsection{Explicit rational block bounds and proof of Lemma \ref{lem:ninecertificate}}\label{subsubsec:nine_bounds}

Let us now check by explicit rational arithmetic that the assumptions of Lemma~\ref{lem:simultaneous} hold on $\mathcal{J}=[177/250,14161/20000]$.

\begin{lemma}[Analytic inertia and positivity bounds]\label{lem:ninecertificate}
Let $a_0=177/250=0.708$, $a_1=14161/20000=0.70805$, $\mathcal{J}=[a_0,a_1]$, $\mu_-:=\mu(a_1)$, and $\mu_+:=\mu(a_0)=150$. Then:
\begin{enumerate}[(i)]
\item For all $a\in\mathcal{J}$, the masses are positive and $\mu'(a)<0$, with $142.08935<\mu_-<142.08937$ (hence $[\mu_-,\mu_+]\subset(142,150]$).
\item The $16\times16$ satellite block $C(a)=\mathscr{Q}(a)|_{\{\xi_0=0\}}$ (and equivalently $\widetilde{K}(a)|_{\{\xi_0=0\}}$) is strictly positive definite for every $a\in\mathcal{J}$.
\item The scalar Schur complement $\sigma(a)$ satisfies $\sigma(a_1)>0$ at $\mu=\mu_-$ and $\sigma(a_0)<0$ at $\mu=\mu_+=150$. Consequently, the endpoint shape Hessians are invertible and have the following Morse indices:
\[
\begin{array}{c|cc}
 \toprule
 \text{Subspace} & \mu=\mu_-\ (a=a_1) & \mu=\mu_+=150\ (a=a_0) \\\midrule
 W_a^{H_{\mathrm{ax}}} & 0 & 1 \\
 W_a^{H_{\mathrm{diag}}} & 0 & 1 \\
 W_a\text{ (full shape space)} & 0 & 2 \\\bottomrule
\end{array}
\]
\end{enumerate}
\end{lemma}
\begin{proof}
Throughout this proof, all terminating decimals denote exact rational numbers. Setting $\beta=\beta_0$ and $w=\omega^2(a)$, we define
$$
 v=R(a)^{-3},\quad \tau_\pm=d_\pm(a)^{-5},\quad
 k=\frac{\sqrt{2}}{2},\quad \delta=\frac{\beta^2}{a^3},\quad
 c_0=S+\beta F(a),\quad c_1=\beta Sv+G(a).
$$
From \eqref{eq:ninemasses}--\eqref{eq:nineomega} we get the relations
\begin{equation}\label{eq:ninecancellation}
 w=\mu+c_0=\mu v+c_1 \quad \textrm{ and } \quad \mu=\frac{c_1-c_0}{1-v},
\end{equation}
which allow us to eliminate $\mu$ from all tangential diagonal entries below. We also set
\begin{align*}
 X&=4\beta\big((a^2-4a+2)\tau_-+(a^2+4a+2)\tau_+\big),\\
 Y&=4\beta\big((a^2+2a-1)\tau_-+(a^2-2a-1)\tau_+\big),\\
 Z&=12a\beta\big((1-a)\tau_--(1+a)\tau_+\big),\\
 h_r&=2(X+Y)-4Z,\qquad h_t=2(X+Y)+4Z.
\end{align*}

\emph{Step 1: Scalar enclosures and monotonicity of $\mu(a)$.}
At $a_0=177/250$ and $a_1=14161/20000$, by squaring the rational numbers $N/10^{12}$ and $(N+1)/10^{12}$ one immediately checks that $N/10^{12}<\zeta<(N+1)/10^{12}$ for each radical $\zeta$ in the table below:
\begin{equation}\label{eq:nineradicalbounds}
\begin{array}{c|r}
 \toprule
 \zeta & N \\\midrule
 \sqrt{2} & 1414213562373 \\
 d_-(a_0) & 765851160474 \\
 d_+(a_0) & 1848926174837 \\
 d_-(a_1) & 765878322581 \\
 d_+(a_1) & 1848991510256 \\\bottomrule
\end{array}
\end{equation}
Inserting the bounds at $a_0$ into \eqref{eq:ninebeta} gives $0.365057017<\beta<0.365057020$. Since $a\mapsto d_\pm(a)=\sqrt{1\pm2a+2a^2}$ is strictly increasing on $\mathcal{J}=[0.708,0.70805]$ (because $2a-1>0$), its values on $\mathcal{J}$ lie between the endpoint bounds in \eqref{eq:nineradicalbounds}; by standard interval arithmetic we thus get, for all $a\in\mathcal{J}$,
\begin{equation}\label{eq:ninecoarse}
\begin{gathered}
 0.365057017<\beta<0.365057020,\qquad 0.9960089<v<0.9962200,\\
 1.62874<c_0<1.62915,\qquad 2.19540<c_1<2.19659,\qquad 140<\mu<151,\\
 -1.48<X<-1.46,\qquad 4.94<Y<4.97,\qquad 3.18<Z<3.21,\\
 0.37<\delta<0.38,\qquad \beta v>0.36.
\end{gathered}
\end{equation}
In particular, $\mu(a)>140>0$ and $w(a)=\mu(a)+c_0(a)>141>0$ on $\mathcal{J}$. Moreover, using $0.7<a<0.71$, $d_->3/4$, $d_+>9/5$, and $v'=-3v/a$, we easily find $|F'(a)|<12$, $|G'(a)|<12$, $-4.3<v'(a)<-4.2$, and $|c_1'(a)-c_0'(a)|<20$ on $\mathcal{J}$, whence
$$
 \mu'(a)=\frac{c_1'(a)-c_0'(a)+\mu(a)v'(a)}{1-v(a)}<\frac{20-140(4.2)}{1-v(a)}<0, \quad \forall\, a\in\mathcal{J},
$$
while evaluating \eqref{eq:ninecancellation} at $a_1$ via \eqref{eq:nineradicalbounds} gives $142.08935<\mu(a_1)<142.08937$. This proves (i).

\emph{Step 2: $D_4$ symmetry blocks of the satellite Hessian $C(a)$.}
We decompose the sixteen-dimensional satellite subspace $\{\xi:\xi_0=0\}\cong\R^{16}$ into $D_4$-isotypic components. Indexing the four vertices of each square by $j=0,1,2,3$ (corresponding to bodies $1,2,3,4$ and $5,6,7,8$), we set
\[
 p=\bigl((1,0),(0,1),(-1,0),(0,-1)\bigr),\qquad
 t=\bigl((1,1),(-1,1),(-1,-1),(1,-1)\bigr),\qquad
 \varepsilon_j=(-1)^j.
\]
Writing $(u,v)$ for satellite displacements on the first and second squares (with $\xi_0=0$), the four two-dimensional character subspaces $V_{\sigma_r,\sigma_s}$ (on which $r$ acts as $\sigma_r\in\{\pm1\}$ and $s$ acts as $\sigma_s\in\{\pm1\}$) have orthogonal bases
\[
\begin{array}{c|cc}
 \toprule
 (\sigma_r,\sigma_s) & \text{First basis vector} & \text{Second basis vector} \\\midrule
 (+,+) & (p,0) & (0,t) \\
 (+,-) & (Jp,0) & (0,Jt) \\
 (-,+) & (\varepsilon p,0) & (0,-\varepsilon Jt) \\
 (-,-) & (-\varepsilon Jp,0) & (0,\varepsilon t) \\\bottomrule
\end{array}
\]
The remaining eight-dimensional orthogonal complement in $\R^{16}$ is the $E$-isotypic component, consisting of four copies of the two-dimensional standard representation $E$. Its $s$-fixed subspace $\R^4$ is spanned by the four basis vectors $e_1,e_2,e_3,e_4$ defined at vertex $j\in\{0,1,2,3\}$ by
\[
\begin{array}{c|cc}
 \toprule
 \text{Basis vector} & \text{First-square displacement at }j & \text{Second-square displacement at }j \\\midrule
 e_1 & \bigl((1+\varepsilon_j)/2,0\bigr) & (0,0) \\
 e_2 & (0,0) & (1,\varepsilon_j) \\
 e_3 & \bigl((1-\varepsilon_j)/2,0\bigr) & (0,0) \\
 e_4 & (0,0) & (1,-\varepsilon_j) \\\bottomrule
\end{array}
\]
while $re_1,re_2,re_3,re_4$ span the $s$-anti-fixed subspace and have the exact same $4\times4$ matrix representation by $r$-invariance of $\mathscr{Q}(a)$. By Schur's lemma, these six subspaces ($V_{++},V_{+-},V_{-+},V_{--}$ of dimension $2$ and the two copies of $\R^4$ from $4E$) are mutually orthogonal with respect to $\mathscr{Q}(a)$ and account for all $2+2+2+2+4+4=16$ satellite dimensions.

Evaluating \eqref{eq:pairhessian} and \eqref{eq:liftedmechanical} on the $s$-fixed standard basis $(e_1,e_2,e_3,e_4)$ yields the $4\times4$ matrix
\[
 Q_E(a)=\begin{pmatrix}R_0&V_0\\V_0^T&T_0\end{pmatrix},
\]
where the $2\times2$ blocks are
\begin{align}
 R_0&=\begin{pmatrix}
 6\mu+2c_0+k+X & -X+Z \\
 -X+Z & 24\beta\mu v+8\beta c_1+h_r+\delta
 \end{pmatrix},\label{eq:nineradial}\\
 V_0&=\begin{pmatrix}-k & -X-Z \\ -Y+Z & -\delta\end{pmatrix},\qquad
 T_0=\begin{pmatrix}
 2c_0+k+Y & -Y-Z \\
 -Y-Z & 8\beta c_1+h_t+\delta
 \end{pmatrix}.\label{eq:ninetangential}
\end{align}
Since $u=MJq$ belongs to the $(+,-)$ character space, the rotational rank-one term in \eqref{eq:liftedmechanical} vanishes on $Q_E$, $C_{++}$, $C_{-+}$, and $C_{--}$. On the four character spaces $V_{\sigma_r,\sigma_s}$, the $2\times2$ matrices of $\mathscr{Q}(a)$ are given by
\begin{align*}
 C_{++}&=\begin{pmatrix}
 12\mu+4c_0+2+4\sqrt{2}+2X & \gamma_0 \\
 \gamma_0 & 24\beta\mu v+8\beta c_1+h_r+(4+\sqrt{2})\delta
 \end{pmatrix},\\
 C_{-+}&=\begin{pmatrix}
 12\mu+4c_0+2-2\sqrt{2}+2X & \gamma_1 \\
 \gamma_1 & 8\beta c_1+h_t+(4-\sqrt{2}/2)\delta
 \end{pmatrix},\\
 C_{--}&=\begin{pmatrix}
 4c_0+2Y-1+4\sqrt{2} & \gamma_2 \\
 \gamma_2 & 24\beta\mu v+8\beta c_1+h_r+(\sqrt{2}-2)\delta
 \end{pmatrix},\\
 C_{+-}&=\ell\begin{pmatrix}1 & -1/a \\ -1/a & 1/a^2\end{pmatrix}
 +\frac{w}{I}\begin{pmatrix}4 \\ 8a\beta\end{pmatrix}\begin{pmatrix}4 & 8a\beta\end{pmatrix},
\end{align*}
where $\ell=2Y+4\beta F>0$ and
\begin{align*}
 \gamma_0&=8\beta\bigl((-4a^2+7a-2)\tau_-+(4a^2+7a+2)\tau_+\bigr),\\
 \gamma_1&=8\beta\bigl((2a^2+a-2)\tau_-+(-2a^2+a+2)\tau_+\bigr),\\
 \gamma_2&=8\beta\bigl((4a^2-a-1)\tau_-+(-4a^2-a+1)\tau_+\bigr).
\end{align*}

\emph{Step 3: Uniform positive definiteness of $C(a)$ on $\mathcal{J}$.}
For $C_{+-}$, both rank-one summands are positive semidefinite and their kernels (spanned by $(1,a)^T$ resp.\ $(2a\beta,-1)^T$) are linearly independent, so $C_{+-}>0$. For $C_{++},C_{-+},C_{--}$, from \eqref{eq:ninecoarse} we immediately get
$$
\begin{array}{c|ccc}
 \toprule
 \text{Block} & \text{First diagonal }> & \text{Second diagonal }> & |\text{Off-diagonal}|< \\\midrule
 C_{++} & 1600 & 1200 & 12 \\
 C_{-+} & 1600 & 26 & 3 \\
 C_{--} & 20 & 1200 & 4 \\\bottomrule
\end{array}
$$
which shows that all three $2\times 2$ blocks have positive trace and determinant on $\mathcal{J}$. As for $Q_E$, the estimates \eqref{eq:ninecoarse} show that the diagonal entries of $R_0$ are larger than $840$ and $1200$ while $|(R_0)_{12}|<5$, whence $R_0>800I_2$; similarly, Gershgorin's theorem applied to $T_0$ gives $T_0>\min(2c_0+k+Y-|Y+Z|,\,8\beta c_1+h_t+\delta-|Y+Z|)I_2>\frac{7}{10}I_2$, while the entries of $V_0$ are bounded in absolute value by $0.71,1.75,1.79,0.38$, so that $\norm{V_0}_{\mathrm{F}}^2<7$. Therefore,
$$
 T_0-V_0^TR_0^{-1}V_0>\left(\frac{7}{10}-\frac{7}{800}\right)I_2>0,
$$
which implies $Q_E(a)>0$ and hence $C(a)>0$, $\forall\,a\in\mathcal{J}$.

\emph{Step 4: Endpoint signs of $\sigma(a)$ and restricted Morse indices.}
In order to find the sign of $\sigma(a)$ without inverting the $16\times16$ matrix $C(a)$ in absolute coordinates, we take advantage of the translational invariance of $U$ along the $s$-fixed subspace of $4E$. More precisely, let $z\in\R$ be the horizontal displacement of the central body and let $x\in\R^4$ denote the absolute satellite coordinates in the basis $(e_1,e_2,e_3,e_4)$. We introduce the relative satellite coordinates $y=x-dz\in\R^4$, where $d=(1,1/2,1,1/2)^T$ (so that $x=y+dz$ shifts every satellite horizontally by $z$). Since $U$ is translation-invariant, $D^2U(q)$ depends only on $y$, whereas the mass term $w\sum_{i=0}^8m_i|\xi_i|^2$ takes the form
$$
 w\big(M_{\mathrm{tot}}z^2+2z\,t_E^Ty+y^TM_Ey\big),\qquad
 M_{\mathrm{tot}}=\mu+4+4\beta,\qquad t_E=(2,4\beta,2,4\beta)^T,
$$
with $M_E=\diag(2,8\beta,2,8\beta)$. Because the rotational lift vanishes on $E$ and $D_4$-equivariance rules out any coupling with the character spaces $V_{\sigma_r,\sigma_s}$, in the coordinates $(z,y)\in\R\oplus\R^4$ the quadratic form is represented by the $5\times5$ matrix
$$
 \widehat Q(a)=\begin{pmatrix}w(a)M_{\mathrm{tot}}&w(a)t_E^T\\w(a)t_E&Q_E(a)\end{pmatrix}.
$$
If we eliminate $z$ in the relative coordinates $(z,y)$, we are left with the $4\times4$ matrix
\begin{equation}\label{eq:nineshapeE}
 Q_{\mathrm{rel}}(a):=Q_E(a)-\theta(a)\,t_Et_E^T,\qquad
 \theta(a)=\frac{w(a)}{\mu(a)+4+4\beta},
\end{equation}
and hence $\det\widehat Q(a)=w(a)M_{\mathrm{tot}}\det Q_{\mathrm{rel}}(a)$. Conversely, since $(z,y)\mapsto(z,x)=(z,y+dz)$ has determinant $1$ and keeps $z$ fixed, eliminating $x$ in the absolute coordinates $(z,x)$ gives $\det\widehat Q(a)=\sigma(a)\det Q_E(a)$, so that
\begin{equation}\label{eq:sigma_det_identity}
 \sigma(a)=w(a)\big(\mu(a)+4+4\beta\big)\Big(1-\theta(a)\,t_E^TQ_E(a)^{-1}t_E\Big)
 =w(a)\big(\mu(a)+4+4\beta\big)\frac{\det Q_{\mathrm{rel}}(a)}{\det Q_E(a)}.
\end{equation}
To determine the sign of $\det Q_{\mathrm{rel}}(a)$, write $e=(2,4\beta)^T\in\R^2$ so that $t_E=(e,e)$, and partition $Q_{\mathrm{rel}}(a)=\bigl(\begin{smallmatrix}R&V\\V^T&T\end{smallmatrix}\bigr)$ into $2\times2$ blocks:
\[
 R=R_0-\theta ee^T,\qquad V=V_0-\theta ee^T,\qquad T=T_0-\theta ee^T,\qquad
 D(a):=T-V^TR^{-1}V.
\]
Since $c_0<4+4\beta$, we have $0<\theta(a)<1$ and $|e|^2=4+16\beta^2<7$, so $R>800I_2-7I_2=793I_2>0$ on $\mathcal{J}$. Hence $\det Q_{\mathrm{rel}}(a)=\det R(a)\cdot\det D(a)$ with $\det R(a)>0$, and the entries of the $2\times2$ Schur complement $D(a)$ are given explicitly by
\begin{equation}\label{eq:nineDentries}
 D_{ij}=T_{ij}-\frac{R_{22}V_{1i}V_{1j}-R_{12}(V_{1i}V_{2j}+V_{2i}V_{1j})+R_{11}V_{2i}V_{2j}}{R_{11}R_{22}-R_{12}^2},\qquad i,j\in\{1,2\}.
\end{equation}

Propagating the rational radical enclosures \eqref{eq:nineradicalbounds} at $a_0=177/250$ and $a_1=14161/20000$ yields the intermediate intervals
\[
\begin{array}{c|cc}
 \toprule
 \text{Quantity} & a=a_0\ (\mu=150) & a=a_1\ (\mu=\mu_-) \\\midrule
 c_0 & (1.62901736,\,1.62901737) & (1.62887046,\,1.62887048) \\
 c_1 & (2.19602524,\,2.19602525) & (2.19596135,\,2.19596136) \\
 X & (-1.47264940,\,-1.47264938) & (-1.47308569,\,-1.47308567) \\
 Y & (4.95445589,\,4.95445594) & (4.95452183,\,4.95452188) \\
 Z & (3.19227954,\,3.19227958) & (3.19134303,\,3.19134307) \\\bottomrule
\end{array}
\]
and the resulting enclosures for the entries of $D(a)$:
\begin{equation}\label{eq:nineDtable}
\begin{array}{c|cc}
 \toprule
 \text{Entry} & a=a_0\ (\mu=150) & a=a_1\ (\mu=\mu_-) \\\midrule
 D_{11} & (4.978,\,4.979) & (4.981,\,4.982) \\
 D_{12} & (-11.028,\,-11.027) & (-11.025,\,-11.024) \\
 D_{22} & (24.414,\,24.415) & (24.410,\,24.411) \\\bottomrule
\end{array}
\end{equation}
From \eqref{eq:nineDtable}, we see that $D_{11}(a_0)>0$ and $D_{11}(a_1)>0$, whereas the determinants at the two endpoints have strictly opposite signs:
\begin{align*}
 \det D(a_0)&=D_{11}(a_0)D_{22}(a_0)-D_{12}(a_0)^2<\frac{4979\cdot24415-11027^2}{10^6}=-\frac{32444}{10^6}<0,\\
 \det D(a_1)&=D_{11}(a_1)D_{22}(a_1)-D_{12}(a_1)^2>\frac{4981\cdot24410-11025^2}{10^6}=\frac{35585}{10^6}>0.
\end{align*}
By \eqref{eq:sigma_det_identity}, $\sigma(a)$ has the same sign as $\det D(a)$, so $\sigma(a_0)<0$ and $\sigma(a_1)>0$. Applying Lemma~\ref{lem:simultaneous}(ii) yields all endpoint Morse indices in part~(iii).
\end{proof}

\begin{figure}[htbp]\centering
\begin{tikzpicture}[>=Latex,font=\small]
\begin{scope}[xshift=-3.5cm,scale=1.3]
\draw (1,0)--(0,1)--(-1,0)--(0,-1)--cycle;
\draw[dashed] (.708,.708)--(-.708,.708)--(-.708,-.708)--(.708,-.708)--cycle;
\foreach \x/\y in {1/0,0/1,-1/0,0/-1} {\fill (\x,\y) circle(1.8pt);}
\foreach \x/\y in {.708/.708,-.708/.708,-.708/-.708,.708/-.708}
 {\draw[fill=white] (\x,\y) circle(2pt);}
\fill (0,0) circle(2.5pt);
\draw[->,thick] (0,0)--(.65,0);
\draw[->,thick,dashed] (0,0)--(.46,.46);
\node[below,align=center] at (0,-1.2)
 {two interlaced squares\\axial and diagonal central displacements};
\end{scope}
\begin{scope}[xshift=-.4cm]
\draw[->] (0,0)--(5.6,0) node[right] {$\mu$};
\draw (.5,.1)--(.5,-.1) node[below] {$\mu_-$};
\draw (4.8,.1)--(4.8,-.1) node[below] {$150$};
\draw[thick] (.5,.7)--(4.8,.7);
\fill (2.6,.7) circle(2.5pt);
\node[above,align=center] at (2.6,.85)
 {common critical mass $\mu_*$\\exact stabilizers $H_{\rm ax}$ and $H_{\rm diag}$};
\node[below,align=center] at (2.6,-.65)
 {restricted indices: $0\to1$ on $W^{H_{\rm ax}},W^{H_{\rm diag}}$\\full shape index: $0\to2$ on $W$};
\end{scope}
\end{tikzpicture}
\caption{Geometry of the nine-body interlaced-square configuration (left, showing the axial and diagonal displacement directions of the central mass) and the Morse index transition across $[\mu_-,150]\subset(142,150]$ (right).}
\label{fig:ninebody}
\end{figure}
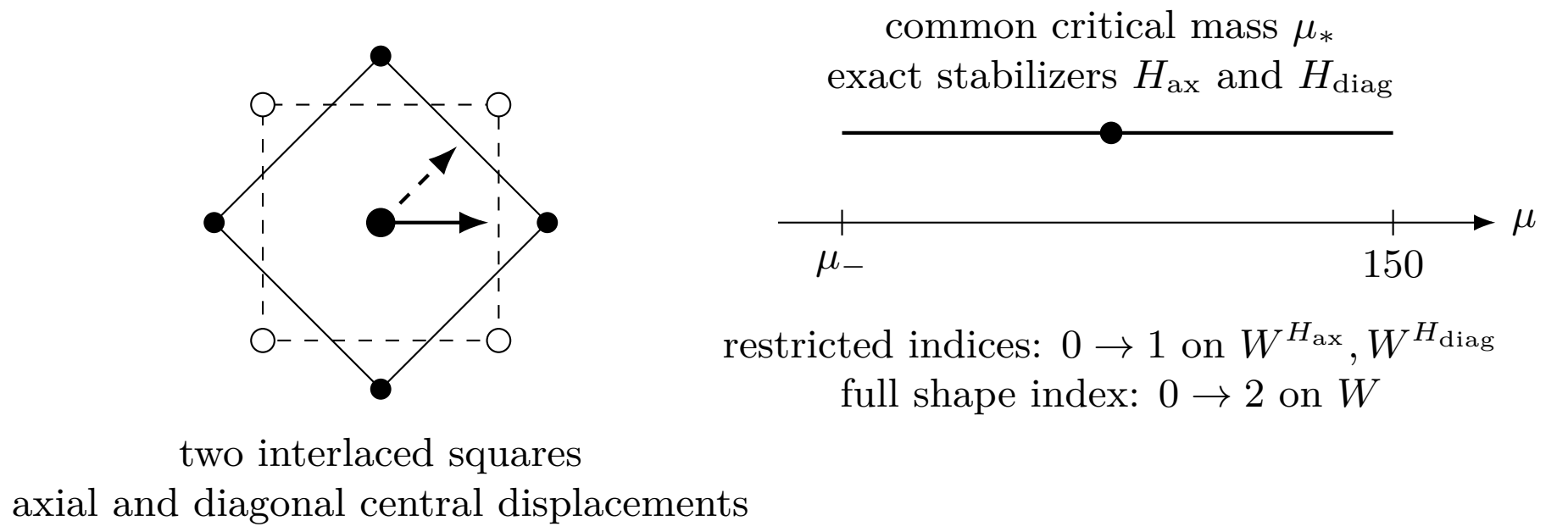

\begin{proof}[Proof of Theorem B]
By Lemma~\ref{lem:ninecertificate}(i), the central mass $\mu(a)$ is strictly decreasing on $\mathcal{J}=[a_0,a_1]$, so we may reparametrize the branch \eqref{eq:ninepositions}--\eqref{eq:ninemasses} by $\mu\in[\mu_-,\mu_+]\subset(142,150]$ and identify the constrained shape spaces $W_a$ with the fixed $D_4$-equivariant space $W:=W_{a_0}$ via the Kato isomorphisms $\Upsilon_a:W\xrightarrow{\sim}W_a$ from Section~\ref{subsec:cc}.
\begin{enumerate}[(i)]
\item By Lemma~\ref{lem:ninecertificate}(iii), as $\mu$ increases from $\mu_-$ to $\mu_+=150$, the restricted Morse indices on both $W^{H_{\mathrm{ax}}}$ and $W^{H_{\mathrm{diag}}}$ increase from $0$ to $1$ (giving spectral flow $s_{H_{\mathrm{ax}}}=s_{H_{\mathrm{diag}}}=0-1=-1$), and the full shape Morse index on $W$ increases from $0$ to $2$ (spectral flow $0-2=-2$).
\item By Lemma~\ref{lem:ninecertificate}(ii)--(iii), the satellite block $C(a)$ is strictly positive definite for all $a\in\mathcal{J}$ and $\sigma(a_0)<0<\sigma(a_1)$. Lemma~\ref{lem:simultaneous}(ii)--(iii) therefore ensures that $\cE_{H_{\mathrm{ax}}}=\cE_{H_{\mathrm{diag}}}=\varnothing$ ($B_{H_{\mathrm{ax}}}=B_{H_{\mathrm{diag}}}=0$), $M_{H_{\mathrm{ax}}}=M_{H_{\mathrm{diag}}}=1$, and there exists a common critical parameter $a_*\in(a_0,a_1)$ (corresponding to a common critical central mass $\mu_*=\mu(a_*)\in(\mu_-,150)\subset(142,150)$) at which $s_{H_{\mathrm{ax}}}(a_*)=s_{H_{\mathrm{diag}}}(a_*)\ne0$. Consequently, both exact-$H_{\mathrm{ax}}$ and exact-$H_{\mathrm{diag}}$ central configurations bifurcate simultaneously from the concentric branch at $\mu_*$ (Figure~\ref{fig:ninebody}).
\item It remains to check that $q_0\ne0$ along every non-trivial bifurcating sequence sufficiently close to the branch. For this purpose, inside $\Sigma_a$ we consider the subset $\Sigma_{a,0}:=\{q\in\Sigma_a:q_0=0\}$. If we take a vector in $\R^{18}$ supported at $q_0$ and project it orthogonally to $\mathcal T_{\mathrm{tr}}$, its $0$-th component is multiplied by $(4+4\beta_0)/(\mu(a)+4+4\beta_0)>0$ while its $M(a)$-scalar products with $\bar q(a)$ and $J\bar q(a)$ are zero; this shows that $\xi\mapsto\xi_0$ maps $W_a$ onto $\R^2$ with $12$-dimensional kernel $W_{a,0}:=W_a\cap\{\xi_0=0\}$. Since $\phi_a(\eta)_0=\eta_0/\sqrt{I_a(\bar q(a)+\eta)}$ in \eqref{eq:radial_chart}, $\phi_a$ maps $B_{r_a}(0)\cap W_{a,0}$ diffeomorphically onto a neighborhood of $\bar q(a)$ in $\Sigma_{a,0}$, and the Hessian of $U_a\circ(\phi_a|_{W_{a,0}})$ at $\eta=0$ is $I(a)^{3/2}C(a)|_{W_{a,0}}>0$ by Lemma~\ref{lem:ninecertificate}(ii). The spaces $W_{a,0}$ form a real-analytic bundle of rank $12$. After pulling back the restricted potentials by a local trivialization near $a_*$, Lemma~\ref{lem:uniqueness} applies on the fixed space $W_{a_*,0}$. Thus, $\bar q(a)$ is locally isolated among critical points of $U_a|_{\Sigma_{a,0}}$ for $a$ near $a_*$. Since any critical point of $U_a|_{\Sigma_a}$ having $q_0=0$ is also critical for $U_a|_{\Sigma_{a,0}}$, we conclude that $q_0\ne0$ on all nearby non-trivial solutions; finally, $H$-invariance forces $q_0\in(\R^2)^H\setminus\{0\}$, i.e.\ along the corresponding reflection axis.\qedhere
\end{enumerate}
\end{proof}

\begin{remark}\label{rem:ninebody_discussion}
In rotating coordinates with angular velocity $\omega$, the linearized equations around a planar central configuration read
$$
 M\ddot{\xi}+2\omega MJ\dot{\xi}=(D^2U(q)+\omega^2M)\xi,
$$
and their restriction to the $J$-invariant shape slice $W$ is governed by the quadratic pencil $P_W(z)=z^2M_W+2\omega z(MJ)_W-K_W$. Since $\dim W=2N-4$ is even and $\det P_W(z)>0$ as $z\to+\infty$, if $K_W$ is non-degenerate and has odd Morse index, then $\det P_W(0)=\det(-K_W)=\det K_W<0$, which implies spectral instability \cite{HuSun2009}. By contrast, an even shape Morse index (such as $m^-(K|_W)=2$ at $\mu=150$) satisfies $\det P_W(0)>0$ and is compatible with gyroscopic stabilization. It is worth noticing that the odd restricted Morse index $m^-(K|_{W^H})=1$ on each reflection-fixed subspace $W^{H_{\mathrm{ax}}}$ and $W^{H_{\mathrm{diag}}}$ does not force dynamical instability, because the symplectic rotation generator $J$ anticommutes with reflections and therefore does not leave $W^H$ invariant.
\end{remark}

\FloatBarrier
\subsection{Periodic and relative-periodic motions: a vertically oscillating ring}\label{subsec:periodic}

To illustrate how the fixed-space method applies to periodic and relative-periodic orbits on an infinite-dimensional loop space, we consider a regular planar ring of $n\ge3$ unit masses surrounding a central body of mass $\mu>0$, where the ring and the central body oscillate in opposite vertical directions in $\R^3$ (compare \cite{GarciaAzpeitiaIze}).

Set
\[
 c_n=\frac{1}{4}\sum_{k=1}^{n-1}\csc\frac{k\pi}{n},\qquad \kappa=\frac{n\mu}{n+\mu}.
\]
With the center of mass fixed at the origin, the ring radius $r>0$, azimuthal angle $\varphi$, and vertical separation $h\in\R$ specify the particle positions
$$
 q_j=\left(r\cos\Big(\varphi+\frac{2\pi j}{n}\Big),\,r\sin\Big(\varphi+\frac{2\pi j}{n}\Big),\,\frac{\mu h}{n+\mu}\right)\quad(1\le j\le n),\qquad
 q_0=\left(0,0,-\frac{nh}{n+\mu}\right).
$$
If we fix the vertical angular momentum $\ell=nr^2\dot{\varphi}$ at $\ell^2=n^2(c_n+\mu)$, Routh's reduction gives the Lagrangian
\begin{equation}\label{eq:ringrouth}
 \mathscr{L}_\ell(r,h,\dot{r},\dot{h})=\frac{n}{2}\dot{r}^2+\frac{\kappa}{2}\dot{h}^2+F(r,h),\qquad
 F(r,h)=\frac{nc_n}{r}+\frac{n\mu}{\sqrt{r^2+h^2}}-\frac{\ell^2}{2nr^2},
\end{equation}
for which $(r,h)=(1,0)$ is an equilibrium corresponding to the planar circular relative equilibrium, with radial and vertical linearized frequencies $\omega^2=c_n+\mu$ resp.\ $\Omega^2=n+\mu$.

\begin{figure}[htbp]\centering
\begin{tikzpicture}[>=Latex]
\draw (-2,0) ellipse(1.3 and .4);
\foreach \t in {0,60,120,180,240,300}
 {\fill ({-2+1.3*cos(\t)},{.4*sin(\t)}) circle(1.7pt);}
\fill (-2,-1.15) circle(2.7pt);
\draw[<->] (-.45,-1.15)--(-.45,0) node[midway,right] {$h(t)$};
\draw[->] (-2,.55)--(-2,1.0);
\draw[->] (-2,-1.35)--(-2,-1.8);
\node[align=left,anchor=west] at (1,-.5)
 {counter-oscillating ring and central body\\
 bifurcating period $T\to 2\pi/\sqrt{n+\mu}$};
\end{tikzpicture}
\caption{The vertically oscillating ring in Proposition~\ref{prop:ringperiodic}, with vertical separation $h(t)$ between the regular $n$-gon and the central mass $\mu$.}
\label{fig:ring}
\end{figure}
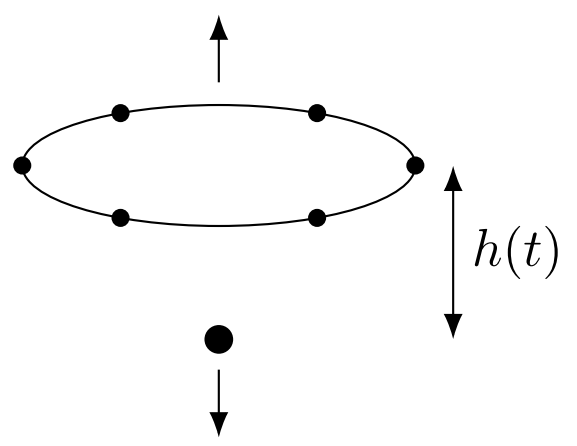

\begin{proposition}[Local oscillating-ring bifurcation]\label{prop:ringperiodic}
Suppose that $\omega/\Omega=\sqrt{(c_n+\mu)/(n+\mu)}\notin 2\mathbb N$. Then a smooth local curve of non-constant, collision-free $2\pi/\nu$-periodic solutions of \eqref{eq:ringrouth} with $h\not\equiv0$ bifurcates from $(\nu,r,h)=(\Omega,1,0)$, with periods tending to $2\pi/\Omega$. In the rescaled time $\tau=\nu t$, such solutions satisfy
$$
 r(-\tau)=r(\tau)=r(\tau+\pi) \quad \textrm{ and } \quad h(-\tau)=h(\tau)=-h(\tau+\pi),
$$
and lift to spatial relative-periodic orbits of the $(n+1)$-body problem (which are periodic in the inertial frame whenever $\Delta\varphi=\int_0^{2\pi/\nu}\frac{\ell}{nr(t)^2}\,dt\in 2\pi\mathbb Q$).
\end{proposition}
\begin{proof}
We work on $\cH=H^1(\R/2\pi\mathbb{Z};\R^2)$ near $(1,0)$ with the action functional
$$
 f_\nu(r,h)=\int_0^{2\pi}\left(\frac{\nu^2}{2}(nr'^2+\kappa h'^2)+F(r,h)\right)d\tau,
$$
and consider the subgroup $H\cong C_2\times C_2$ generated by $\tau\mapsto-\tau$ and $(\tau,r,h)\mapsto(\tau+\pi,r,-h)$. Elements in $\cH^H$ are of the form $r(\tau)=1+\sum_{k\ge0}a_k\cos(2k\tau)$ and $h(\tau)=\sum_{k\ge0}b_k\cos((2k+1)\tau)$. Since $F_{rr}(1,0)=-n\omega^2$, $F_{hh}(1,0)=-\kappa\Omega^2$, and $F_{rh}(1,0)=0$, the quadratic form $D^2f_\nu(1,0)$ on $\cH^H$ is diagonal in this Fourier basis with coefficients $n((2k\nu)^2-\omega^2)$ and $\kappa(((2k+1)\nu)^2-\Omega^2)$, $k\ge0$. Since $\omega/\Omega\notin2\mathbb{N}$, at $\nu=\Omega$ the restricted kernel on $\cH^H$ is one-dimensional and spanned by $(0,\cos\tau)$, with positive crossing derivative $2\kappa\Omega\pi>0$. Applying Corollary~\ref{cor:simple} on $\cH^H$ and Theorem~\ref{thm:PSC} yields a smooth local curve of critical points with $h(s)=s\cos\tau+o(s)\not\equiv0$ for small $s\ne0$ (Figure~\ref{fig:ring}). Since the cyclic ring configuration space is invariant under the full Newtonian equations, Routh reconstruction via $\dot\varphi=\ell/(nr^2)$ lifts these solutions to spatial relative-periodic orbits. Notice that for $n=8$ equal satellites one has $\sin(k\pi/8)\ge\sin(\pi/8)>1/4$, whence $c_8<7<8$ and $\omega/\Omega<1$, $\forall\,\mu>0$.
\end{proof}

It is worth observing that in Proposition~\ref{prop:ringperiodic} the restriction to the discrete fixed-point space $\cH^H$ automatically rules out the zero modes due to time shifts on $H^1(\R/2\pi\mathbb{Z};\R^2)$ (cf.\ \cite{ChencinerMontgomery,FukudaOzaki,HuPortaluriYang,HuSun}).

\section{Concluding remarks and open problems}\label{sec:conclusion}

We finish the paper by mentioning two open questions that naturally arise from our study of the nine-body family \eqref{eq:ninepositions}--\eqref{eq:ninemasses}.

The first question concerns the uniqueness and transversality of the critical mass $\mu_*\in(\mu_-,150)$. While Lemma~\ref{lem:simultaneous} and Theorem B show that any zero of $\sigma(a)$ with non-zero local spectral flow is a simultaneous bifurcation instant for both $H_{\mathrm{ax}}$ and $H_{\mathrm{diag}}$, our endpoint estimates do not yet prove that $\sigma'(a)>0$ across the whole interval $\mathcal J$. On the $2$-dimensional kernel $E_{a_*}=\ker(K(a_*)|_{W_{a_*}})\cong E$, by Theorem~\ref{thm:LS} the $D_4$-invariant Lyapunov--Schmidt reduced functional $r_a(z_1,z_2)$ around $(a_*,0)$ takes the form
$$
 r_a(z_1,z_2)=\frac{1}{2}\rho(a)\big(z_1^2+z_2^2\big)
  +\frac{\tau_1(a)}{4}\big(z_1^4+z_2^4\big)
  +\frac{\tau_2(a)}{2}z_1^2z_2^2+O(|z|^6),
$$
with $\rho(a_*)=0$ and $\sgn\rho(a)=\sgn\sigma(a)$ for $|a-a_*|<\delta_0$. Comparing the determinant of $D^2(U_a\circ\phi_a)(0)|_{W_a^{H_{\mathrm{ax}}}}$ in the Lyapunov--Schmidt splitting with \eqref{eq:sigma_det_identity} and \eqref{eq:K_def}, we see that $\rho(a)=c(a)\sigma(a)$ for a positive real-analytic factor $c(a)>0$, so that $\rho'(a_*)=c(a_*)\sigma'(a_*)$. Therefore, proving that $\sigma'(a)>0$ on $\mathcal{J}$ and determining the signs of $\tau_1(a_*)$ and $\tau_1(a_*)+\tau_2(a_*)$ would show that $\mu_*$ is unique and establish whether the axial ($z_2=0$) and diagonal ($z_1=\pm z_2$) branches bifurcate supercritically or subcritically with respect to $\mu$.

The second question is related to the linear stability of the relative equilibrium at $\mu(a_0)=150$. On the $14$-dimensional shape slice $W_a$, the gyroscopic pencil $P_W(z)=z^2M_W+2\omega(a)z(MJ)_W-K_W(a)$ satisfies $P_W(z)^T=P_W(-z)$, so that setting $z=i\nu$ gives a real polynomial of degree $14$ in $x=\nu^2$, namely
$$
 \Pi_a(\nu^2):=\det P_W(i\nu)
 =\det\big(-\nu^2M_W+2i\omega(a)\nu(MJ)_W-K_W(a)\big).
$$
Since simple eigenvalues of $P_W(z)$ are automatically semisimple, verifying via rigorous interval arithmetic or analytic bounds that $\Pi_{a_0}(x)=0$ admits $14$ distinct positive roots $0<x_1<\dots<x_{14}$ in $(0,\infty)$ would yield $28$ simple purely imaginary eigenvalues $z=\pm i\sqrt{x_k}\in i\R\setminus\{0\}$ and thereby provide a complete computer-assisted proof of linear stability for the reduced shape dynamics at the $\mu(a_0)=150$ configuration of \cite{Seamons2026}.

\section*{Acknowledgements}

Alessandro Portaluri acknowledges the support of the Istituto
Nazionale di Alta Matematica (INdAM) and its Gruppo Nazionale
per l'Analisi Matematica, la Probabilit\`a e le loro Applicazioni
(GNAMPA). His research was partially supported by GNAMPA--INdAM,
Italy.

\phantomsection
\addcontentsline{toc}{section}{References}
\setlength{\bibsep}{2pt}

\bigskip
\noindent\begin{minipage}{\textwidth}
{\small
\raggedright

\noindent\textbf{Daniele Garrisi}\\
School of Mathematical Sciences, University of Nottingham Ningbo China\\
199 Taikang East Road, Ningbo 315100, China\\
Email: \texttt{daniele.garrisi@nottingham.edu.cn}

\medskip

\noindent\textbf{Alessandro Portaluri}\\
DISAFA, Universit\`a degli Studi di Torino\\
Largo Paolo Braccini 2, 10095 Grugliasco (Torino), Italy\\
Email: \texttt{alessandro.portaluri@unito.it}

\medskip

\noindent\textbf{Li Wu}\\
School of Mathematics, Shandong University\\
Jinan 250100, China\\
Email: \texttt{vvvli@sdu.edu.cn}

\par}
\end{minipage}
\end{document}